\documentclass[a4paper, 12pt]{amsart}
\usepackage{graphics}
\usepackage{ifthen}
\usepackage{amsmath,amsfonts,amstext,amscd,bezier,amssymb,a4,enumerate,epsfig,psfrag,subfigure}
\usepackage{graphicx}
\usepackage{amsthm,amsopn}
\usepackage{bm}
\usepackage{amssymb}
\usepackage[T1]{fontenc}
\usepackage{mathrsfs}

\usepackage{color}

\everymath={\displaystyle}
\let\cal=\mathcal

\newcommand{\e}{{\ = \ }}
\newcommand{\ms}{{\medskip}}
\newcommand{\bs}{{\bigskip}}
\newcommand{\noi}{{\noindent}}
\newcommand{\beql}{\begin{equation}\label}
\newcommand{\beqn}{\begin{equation*}}
\newcommand{\eeq}{\end{equation}}
\newcommand{\eeqn}{\end{equation*}}

\newcommand{\NN}{{\mathbb{N}}}
\newcommand{\ZZ}{{\mathbb{Z}}}

\newcommand{\RR}{{\mathbb{R}}}

\newcommand\diag{\operatorname{diag}}

\newcommand\sn{\operatorname{sn}}
\newcommand\dn{\operatorname{dn}}
\newcommand\spec{\operatorname{spec}}
\newcommand\tr{\operatorname{tr}}

\newtheorem*{theo*} {Theorem}

\newcommand\blfootnote[1]{%
  \begingroup
  \renewcommand\thefootnote{}\footnote{#1}%
  \addtocounter{footnote}{-1}%
  \endgroup
}

\title{On the open Toda chain with external forcing}

\author[P.~Deift, L.-C.~Li, H.~Spohn, C.~Tomei, T.~Trogdon] {Percy Deift, Luen-Chau Li, Herbert Spohn,\\ Carlos Tomei, Thomas Trogdon}
\address{Percy Deift\\
               Courant Institute\\
               251 Mercer St.\\
               New York, NY 10012, USA}
\email{deift@cims.nyu.edu}
\address{Luen-Chau Li\\
               Department of Mathematics\\
               Pennsylvania State University\\
               University Park, PA 16802, USA}
\email{luenli@math.psu.edu}
\address{Herbert Spohn\\
                Zentrum Mathematik\\
                Technische Universit\"at M\"unchen\\
                Boltzmannstrasse 3\\
                D-85748, Garching, Germany}
\email{spohn@ma.tum.de}
\address{Carlos Tomei\\
                Depto. de Matem\'atica\\
                PUC-Rio\\
                Rio de Janeiro, Brazil}
\email{carlos.tomei@gmail.com}
\address{Thomas Trogdon\\
               Department of Applied Mathematics\\
               University of Washington\\
               Seattle, WA 98195, USA}
\email{trogdon@uw.edu}


\dedicatory{In honor of Shmuel Agmon and his many contributions to mathematics}

\date{}

\begin{document}
\maketitle

\centerline{The Duck Test: `If it looks like a duck, walks like a duck and}
\centerline{quacks like a duck... it's a duck!'}

\begin{abstract}
We consider the open Toda chain with external forcing, and in the case 
when the forcing stretches the system, we derive the longtime behavior
of solutions of the chain. Using an observation of J\"{u}rgen Moser, we then show that the system  is completely integrable, in the sense that the $2N$-dimensional system has $N$ functionally independent Poisson commuting integrals, and also has a Lax-Pair formulation. In addition, we construct action-angle variables for the flow. In the case when the forcing compresses the system, the analysis of the flow remains open.
\end{abstract}

\medbreak

%{\noi\bf Keywords:}  Toda lattice, Liouville integrability, wave map, action-angle variables

%\section{Introduction}\label{sec:intro}
In 1967, Morikazu Toda introduced \cite{To1967} the eponymous Toda system with Hamiltonian

\begin{equation}\label{1.1}
H(q,p) \e \frac{1}{2m} \sum_{n \in \ZZ} \ p_n^2 \ + \frac{a}{b} \sum_{n \in \ZZ}  \ e^{-b(q_{n+1} - q_n)} \ + \ a \sum_{n \in \ZZ} \ (q_{n+1} - q_n ),\, \hbox{for \ } a, b >0
\end{equation}
for particles of equal mass $m>0$ with positions $q= \{q_n\}$  on the line, and momenta $p = \{p_n\}$. The Hamiltonian equations generated by $H$ have the form

\begin{equation}\label{1.2}
\begin{array}{ll}
    \dot{q}_n \ &= \  \frac{\partial H}{\partial p_n} \e \frac{1}{m} p_n,  \\
    \dot{p}_n \ &= \ - \ \frac{\partial H}{\partial q_n} \e - a \ \big( e^{-b(q_{n+1} - q_n)} - e^{-b(q_{n} - q_{n-1})}\big),
\end{array}
\end{equation}
for $ - \infty < n < \infty$, and so
\begin{equation}\label{2.1}
\ddot{q}_n \e   - \frac{a}{m} \ \big( e^{-b(q_{n+1} - q_n)} - e^{-b(q_{n} - q_{n-1})}\big) \ .
\end{equation}

\ms
Toda's goal in considering $H$ was to investigate further the observation of J.~ Ford and J.~Waters \cite{FoWa1963}, found by numerical computation, that nonlinear systems have `normal modes where a normal mode is defined as a motion in which each oscillator moves at essentially constant amplitude (energy) and at a given frequency'. As Toda notes, the existence of such normal modes implies, in particular, the non-ergodic character of the system. Toda's specific goal in \cite{To1967} was to discover an exact, explicit form for such normal modes for $H$, and he found (in the case $m=1$) such modes in the form of a travelling wave solution for (\ref{2.1}),
\begin{equation*}
q_{n+1} - q_n \e - \frac{1}{b} \log \left( 1 + \frac{4(K \nu)^2}{ab} \left[ \mathrm{d n}^2 ( 2 K (\nu t \pm \frac{n}{\lambda})) - \frac{E}{K}\right] \right) \, , -\infty < n < \infty \ ,
\end{equation*}
where for given any wave length $\lambda$ and modulus $0 < k < 1$, $K = K(k)$ and $E = E(k)$ are complete elliptic integrals
\begin{equation*}
K \e \int_0^{\pi/2} \ \frac {d \theta}{\sqrt{ 1 - k^2 \sin^2 \theta}} \ , \quad
E \e \int_0^{\pi/2} \sqrt{1 - k^2 \sin^2 \theta} \ d  \theta,
\end{equation*}
and the frequency $\nu$ is given by
\begin{equation*} 
\nu \e \frac{1}{2K} \ \sqrt{ \frac{ab}{ \big( \frac{1}{\sn^2( 2 K/\lambda)} - 1 + \frac{E}{K} \big) }} \ .
\end{equation*}
Here $\sn$ and $\dn$ are the standard Jacobi elliptic functions.

\ms
In the case that  $n$ runs over $\ZZ$ in (\ref{1.1}), the linear term $a \sum_{n \in \ZZ} \ (q_{n+1}- q_n)$ plays no role. Scaling
\[ q_n \to b q_n \ , \quad t \to t \sqrt{ \frac{m}{ab}} \]
in (\ref{2.1}), we see that we can restrict our attention to the case

\begin{equation}\label{3.2}
H \e \frac{1}{2} \sum_{n \in \ZZ} \ p_n^2 \ + \ \sum_{n \in \ZZ} e^{q_n - q_{n+1}} \ + \ c \sum_{n \in \ZZ} \ (q_n - q_{n+1}) \ ,
\end{equation}
where $c$ is any constant.

\ms
In the periodic case,
\[ q_{n+N} \e q_n +s\ , \quad  p_{n+N} \e p_n \]
for some $N \in \NN$ and $s \in \RR$, equations (\ref{1.2}) scaled as above take the form
\begin{equation}\label{3.3}
    \begin{array}{ll}
    \dot{q}_n \ &= \   p_n  \\
    \dot{p}_n \ &= \  \big( e^{(q_{n+1} - q_n)} - e^{(q_{n} - q_{n-1})}\big)
\end{array}
\end{equation}
for $1\le n \le N$, where $q_{N+1} \e q_1 +s $  as above and 
$-s =  \sum_{n=1}^N (q_n - q_{n+1})$. In an extensive numerical investigation in the cases $N=3$ and $N=6$, J.~Ford, S.~ Stoddard and J.~ Turner \cite{FoStTu1973} found strong evidence that the lattice was integrable. And indeed, inspired by \cite{FoStTu1973}, M.~H\'enon \cite{He1974} and, independently, H.~Flaschka \cite{Fl1974} showed that there are $N$ independent integrals for (\ref{3.3}). Also independently, S.~Manakov \cite{Ma1974}\footnote{For the record: the articles by H\'enon and Flaschka were submitted on August 13, 1973 and August 22, 1973, respectively, and the article of Manakov was submitted on February 8, 1974. Also, Henon, Flaschka and Manakov only considered the case with $s= 0$, but their methods go through for general $s$. } proved the same result, viz., there are $N$ independent integrals for (\ref{3.3}). The starting point of H\'enon's analysis was the integrability of the hard-sphere gas, which is a limiting form of the Toda lattice. On the other hand, Flaschka and Manakov based their analysis on the observation that (\ref{3.3}) can be written in Lax-pair form as follows: set
\begin{equation*}
a_i \e - p_i/2 \ , \quad b_i \e \frac{1}{2} e^{(q_i - q_{i+1})/2} \ , 1 \le i \le N \ , \end{equation*}
where
$q_{N+1} \e q_1 +s  ,\  p_{N+1} \e p_1$,
and in the variables $a_i, b_i$, equations (\ref{3.3}) take the form 
\begin{equation}\label{5.1}
\dot{a}_i \e 2 (b_i^2 - b_{i-1}^2) \ , \quad \dot{b}_i \e b_i (a_{i+1} - a_i) \ , \quad 1 \le i \le N \ ,
\end{equation}
where $b_0 = b_N$ and $a_{N+1} = a_1$. Note that the {\it stretch parameter} $s$ is now part of the initial conditions.

Define the symmetric matrix
\begin{equation*}
L \e \begin{pmatrix} a_1 & b_1 & 0 & b_N \\ b_1 & a_2 & \ddots & 0 \\
0 & \ddots & \ddots & b_{N-1} \\ b_N & 0 & b_{N-1} & a_N \end{pmatrix} \ \e \ L^T
\end{equation*}
and the skew-symmetric matrix
\begin{equation*}
B \e \begin{pmatrix} 0 & -b_1 & 0 & b_N \\ b_1 & 0 & \ddots & 0 \\
0 & \ddots & \ddots & -b_{N-1} \\ -b_N & 0 & b_{N-1} & 0 \end{pmatrix} \
 \e - B^T \ .
\end{equation*}
Then, if $q_i(t), p_i(t)$  solve (\ref{3.3}), so that $a_i(t), b_i(t)$ solve (\ref{5.1}), then $L=L(t)$ solves

\begin{equation*}
\dot{L} \e [ L, B ] \e LB - BL
\end{equation*}
with $L_0 \e L(t=0)$ given by $q_i(0), p_i(0)$. By the general theorem of Lax \cite{La1968}, it follows immediately that $t \mapsto L(t)$ is an isospectral deformation, i.e., \[ \spec L(t) \e \spec L_0,\,\, \, t>0 . \]
In particular, the eigenvalues $\lambda_1, \ldots, \lambda_N$ of $L_0$ are $N$ constants of the motion for (\ref{3.3}), and in \cite{Fl1974}, Flaschka relates the $\lambda_i$'s to the integrals of the motion obtained by H\'enon. Furthermore, in \cite{Ma1974}, Manakov showed that the $\lambda_i$'s Poisson commute, so that the Hamiltonian system (\ref{3.3}) is completely integrable in the sense of Liouville. In principle, this meant that the periodic Toda system could be solved up to quadrature, and indeed in \cite{KaMo1975} M. Kac and P. van Moerbeke showed how to use the Lax-pair formalism to partially solve (\ref{3.3})  in terms of hyperelliptic function theory: A full solution was given shortly thereafter by E.Date and S.Tanaka \cite{DaTa1976}.

\ms
The methods of H\'enon, Flaschka and Manakov can also be used to prove the integrability of the Toda lattice with other boundary conditions, particularly 
the so-called open Toda lattice (sometimes referred to as Toda with ``fixed-ends''---see below), and also scattering systems with infinitely many particles (see \cite{He1974}, \cite{Fl1974}, \cite{Ma1974} and also \cite{Fla1974}).

In this paper we are particularly interested in the {\it open Toda lattice} where  the Hamiltonian has the form
\begin{equation*}
H_F(q,p) \e \frac{1}{2} \sum_{n=1}^N \ p_n^2 \ +  \sum_{n =1}^{N-1}  \ e^{(q_{n} - q_{n+1})}
\end{equation*}
giving rise to the equations
\begin{equation}\label{7.2}
\begin{array}{llll}
    \dot{q}_n \ &= \   p_n , \quad 1 \le n \le N,  \\
    \dot{p}_1 \ &= \ - e^{q_1 - q_2} \\
    \dot{p}_n \ &= \  e^{q_{n-1} - q_{n}} - e^{q_{n} - q_{n+1}} , 2 \le n \le N-1, \\
    \dot{p}_N \ &= \ e^{q_{N-1}-q_N} .
\end{array}
\end{equation}
One can view (\ref{7.2}) as arising from (\ref{3.3}) with $N+2$ particles $q_0, q_1, \ldots, q_N, q_{N+1}$ by setting
\[ q_{N+1} \e \infty \ , \quad q_0 = - \infty \]
so that the ends are ``fixed'' at $\pm \infty$: For this reason the open Toda lattice is sometimes referred to as Toda with ``fixed-ends''. Equations (\ref{7.2}) can be written in Lax-pair form by setting
\beql{7.5}
\begin{array}{ll}
a_i \e - p_n/2 ,\,  1 \le n \le N,\\
b_i \e \frac{1}{2} e^{(q_n - q_{n+1})/2} ,\, 1 \le n \le N-1,
\end{array}
\eeq
and defining the symmetric matrix
\begin{equation*}
L_F \e \begin{pmatrix} a_1 & b_1 & \ldots & 0 \\ b_1 & a_2 & \ddots & 0 \\
\vdots & \ddots & \ddots & b_{N-1} \\0 & 0 & b_{N-1} & a_N \end{pmatrix} \ \e \ L_F^T
\end{equation*}
and the skew-symmetric matrix
\begin{equation*}
B_F \e \begin{pmatrix} 0 & -b_1 & \ldots & 0 \\ b_1 & 0 & \ddots & 0 \\
\vdots & \ddots & \ddots & -b_{N-1} \\0 & 0 & b_{N-1} & 0 \end{pmatrix}
 \e - B_F^T \ .
\end{equation*}
Then if $q_n(t), p_n(t)$ solve (\ref{7.2}), $L_F(t)$ solves the Lax-pair equation
\begin{equation}\label{8.3}
\dot{L}_F \e [ L_F, B_F]
\end{equation}
with $L_{F,0} \e L_F(t=0)$ given by $q_n(0), p_n(0)$. Again, the eigenvalues $\lambda_i(t)$ of $L_F(t)$ are constant and hence give $N$ integrals for the open Toda lattice. In what follows, we will often simply refer to (\ref{7.2}) as the   Toda lattice, or the Toda system, or the Toda chain,  and in cases where the periodic problem is under discussion, we will specifically refer to  (\ref{3.3})  as the periodic Toda system.

\ms
Just as in the periodic case, the Toda system can be solved explicitly, now in terms of rational functions of exponentials, as shown by J.~Moser in \cite{Mo1975}. Furthermore, Moser showed that the system has the following remarkable long-term scattering behavior:
\begin{equation}\label{9.1}
q_n(t) \e \alpha_n^\pm t + \beta_n^\pm + O(e^{-\delta|t|}) , \ t \to \pm \infty,\,\,\, \delta > 0 , 1 \le n \le N
\end{equation}
with
\begin{equation}\label{9.2}
\alpha_n^+ \e - 2 \lambda_n, \quad \alpha_n^- \e - 2 \lambda_{N-n+1},  \quad 1 \le n \le N
\end{equation}
and scattering shift as $t$ goes from $-\infty$ to $\infty,$ given by
\begin{equation}\label{9.3}
\beta_{N-n+1}^+ - \beta_n^- \e \sum_{j \ne k} \ \ln (\alpha_j^- - \alpha_k^+)^2 \ .
\end{equation}
Explicit expressions for the $\beta_n^\pm$'s themselves were derived later in the early 2000's (see \cite{DeDuTr2021}), and subsequently in \cite{LeSaTo2008} ; see (\ref{35.3}) below.
Here $\lambda_1 > \lambda_2 > \ldots > \lambda_N$ are the eigenvalues of $L_{F,0}$.

When one of the authors (PD) came across formula (\ref{9.3}), he was astounded: he had just completed a PhD in abstract scattering theory in Hilbert space, and the idea that one could compute the scattering shifts (equivalently, the scattering matrix) for an $N$ particle system explicitly, was beyond anything he had ever encountered. When he asked Moser how this was possible, Moser replied, somewhat mysteriously, that `Every scattering system is integrable!' 

What Moser meant was the following (see \cite{MoZe2005}, Integrals via Asymptotics: the St{\"o}rmer Problem): suppose one has the solution of a Hamiltonian system
\[ (q(t), p(t)) \e (q_1(t), \ldots, q_N(t), p_1(t), \ldots, p_N(t)) \in \RR^{2N}\]
with Hamiltonian $H$ and with the property that, as $t \to \infty$,
\begin{equation*}
\begin{array}{ll}
p(t) \e p_\infty + o(1/t), \\
q(t) \e q_\infty + t p_\infty + o(1),
\end{array}
\end{equation*}
for some constants $(q_\infty, p_\infty)$. Let $U_t(q(0),p(0)) \e (q(t), p(t))$ be the solution of the system with initial data $(q(0),p(0))$ and let
$ U_t^0 (q^0(0), p^0(0)) \e (q^0(t), p^0(t))$, where $(q^0(t), p^0(t))$ solves the free particle motion with Hamiltonian $H^0(q,p) \e p^2/2$, so
\begin{equation*}
\begin{array}{ll}
p^0(t) \e p^0(0) \\
q^0(t) \e q^0(0) + p^0(0) t \ .
\end{array}
\end{equation*}
Then
\begin{equation*}
\begin{array}{lll}
 U_{-t}^0 \circ U_t (q_0,p_0) &\e U_{-t}^0 ( q_\infty + p_\infty t+ o(1), p_\infty + o(1/t)) \\
 &\e (q_\infty + p_\infty t + o(1) - (p_\infty + o(1/t))t, p_\infty + o(1/t)) \\
 &\e (q_\infty + o(1), p_\infty + o(1/t)) \to (q_\infty, p_\infty) \ \quad \hbox{as} \ \ t \to \infty \ .
 \end{array}
 \end{equation*}
Thus the wave operator
\[ W(q_0,p_0) \equiv \lim_{t \to \infty} \ U_{-t}^0 \circ U_t (q_0,p_0) \e (q_\infty, p_\infty)\]
exists. But then
\[ U_{-t}^0 \circ U_t \circ U_s \e U^0_s \circ U_{-(t+s)}^0 \circ U_{t+s} \]
implies
\[ W \circ U_s \e U_s^0 \circ W \]
or,  if $W^{-1}$ exists,
\begin{equation}\label{11.2}
U_s \e W^{-1} \circ U_s^0 \circ W \ .
\end{equation}

Now $U_{-t}^0 \circ U_t$ is symplectic for all $t$ and so $W$, and hence $W^{-1}$, are symplectic. Thus (\ref{11.2}) shows us that $U_s$ is symplectically equivalent to $U_s^0$, and hence is completely integrable. Indeed, if $\alpha_1, \ldots, \alpha_N,$ are commuting integrals for $H^0$, then $\beta_i = \alpha_i \circ W, i=1, \ldots,N$ are commuting integrals for $H$:
\[ \beta_i \circ U_t(q(0),p(0)) \e \alpha_i \circ W \circ U_t(q(0),p(0)) = \alpha_i \circ U_t^0(W(q(0),p(0)) \e \hbox{constant} \]
and as $W$ is symplectic,
\begin{equation}\label{12.1}
\{ \beta_i, \beta_j \} \e \{ \alpha_i \circ W , \alpha_j \circ W \} \e \{ \alpha_i, \alpha_j \} \circ W \e 0 \ .
\end{equation}
Said differently, the above calculation shows more generally that `if a system behaves like an integrable system, then it is an integrable system!' or, as in the famous `Duck Test', `if it looks like a duck, walks like a duck and quacks like a duck... it's a duck!'

In Moser's argument in \cite{Mo1975} one finds that in addition to (\ref{9.1}) one also has
\begin{equation} \label{13.1}
p_n \e \alpha_n^\pm + O(e^{-\delta |t|}) \ \quad \hbox{as} \ \ t \to \pm \infty
\end{equation}
and so by Moser's integrability argument, the Toda lattice is integrable. Note further that
\[ \alpha_i(q,p) \e p_i \ , 1 \le i \le N \]
are commuting integrals of the motion for $H^0$ and so
\[ \beta_i(q_0,p_0) = \alpha_i(W(q_0,p_0)) \e \alpha_i(q_\infty, p_\infty) \e p_{\infty,i},\,\, 1\leq i\leq N \]
are commuting integrals of the motion for $H_F$. But, from (\ref{9.1}) and (\ref{9.2}), up to a factor of $-2$, the $p_{\infty,i}$'s are just the eigenvalues of $L_{F,0}$, as they should be!

From the ``duck'' we learn that there is an interesting Catch 22 in the problem: We could not have derived, by any means, utilizing any and all dynamical tools, the asymptotic behavior of the system, unless it was integrable in the first place!

Moser's argument can be used to prove the integrability of a variety of dynamical systems. For example, in \cite{DeZh2002} the authors showed, contrary to expectations, that the perturbed defocusing nonlinear Schr\"odinger equation on the line,
\begin{equation*}
\begin{array}{ll}
i q_t + q_{xx} - 2 |q|^2 q - \epsilon K(|q|^2)q = 0 \\
q(x,t=0) = q_0(x) \to 0 \ , \quad \hbox{as} \ \ |x| \to \infty
\end{array}\end{equation*}
is integrable for $0 < \epsilon< \epsilon_0$, for some $\epsilon_0 > 0$. Here $K(|q|^2) = O(|q|^{\ell})$ as $|q| \to 0$ for suitably large $\ell >2$.
 
Wave operators were first introduced in the context of quantum mechanical scattering theory, by C.M\o{}ller \cite{Mo1945} . In particular, for two self-adjoint operators
$A$ and $B$ in Hilbert space, with propagators $e^{iAt}$ and $e^{iBt}$, M\o{}ller introduced the quantum mechanical wave operator
$$ Wf \e \lim_{t \to \infty} e^{-iAt}e^{iBt} f$$
for vectors $f$ in the Hilbert space. Cook \cite{Co1967}, and shortly thereafter Hunziker \cite{Hu1968}, were the first  to use wave operators in the context of classical dynamics. 
Wave operators can also be used more generally to address problems in analysis:  see, in particular, Nelson's proof in \cite{Ne1969} of Sternberg's Linearization Theorem for non-resonant systems, and also the proof of Darboux's Theorem establishing the existence of local canonical coordinates for symplectic forms in \cite{MoZe2005}. A detailed treatment of classical two-body scattering theory in three dimensions using wave operators is given in \cite{Si1971}. A detailed treatment of classical and quantum mechanical 
N-particle scattering theory is given in \cite{DeGe1997}.

\bs
In this paper we consider Toda's original system (\ref{3.2}) in the finite fixed-end case, with $c \ne 0$. The Hamiltonian for the system has the form
\begin{equation*}
H_c(q,p) \e \frac{1}{2} \sum_{n=1}^N p_n^2 + \sum_{n=1}^{N-1} \ e^{q_n - q_{n+1}} + c \ \sum_{n=1}^{N-1} (q_n - q_{n+1}),
\end{equation*}
giving rise to the associated Hamiltonian equations
\begin{equation}\label{14.2}
\begin{array}{llll}
\dot{q}_n &\e p_n \ , \quad 1 \le n \le N, \\
\dot{p}_1 &\e - e^{q_1 - q_2} - c, \\
\dot{p}_n &\e e^{q_{n-1} - q_n} - e^{q_n - q_{n+1}} \ , \quad 2 \le n \le N-1, \\
\dot{p}_N &\e e^{q_{N-1} - q_N} + c .
\end{array}
\end{equation}
As opposed to the whole line case and the periodic case, the constant $c$, and in particular the sign of $c$, now plays a determining role. Note that
\[ c \sum_{n=1}^{N-1} (q_n - q_{n+1}) \e c (q_1 - q_N) , \]
and we think of $H_c$ as the Hamiltonian of a lattice of particles $q_1, \ldots, q_N$ with external forces
acting on the endpoints of the lattice via the potential $cq_1 - cq_N$. When $c>0$, the forces
\[ - \frac{\partial}{\partial q_1} c (q_1 - q_N) \e - c \ \ \quad - \frac{\partial}{\partial q_N} c (q_1 - q_N) = c \]
stretch the lattice , and when $c<0$, they compress the lattice.

The study of $H_c$ is motivated
in part by the statistical mechanics of the Toda lattice. Here the statistics of the statistical mechanical ensemble is given by the canonical  measure
$$  \prod_{j=1}^N \mathrm{d}p_j \prod_{j=1}^{N} \mathrm{d}q_j \exp ( -\beta H_c(q,p))$$
suitably normalized: This measure is clearly invariant under the flow generated by $H_c$.  From \cite{To1989} we learn that the case $c<0$ arises most naturally, but the case $c>0$ is also of interest. In general, the analysis of the thermodynamic limit, $N \to \infty$, of a statistical 
mechanical system with Hamiltonian $H$, is greatly simplified if $H$
is known to be integrable. This motivates, in particular, the study of the integrability of $H_c$.

The numerical calculations below suggest strongly that in the case  $c>0$, $H_c$ is integrable. And indeed, the main result in this paper is to show, using Moser's integrability argument, that this is the case. In the case $c <0$, we will argue below that the numerical calculations suggest that also in this case there is integrable structure, or near integrable structure, associated with the system.

\begin{figure}[tbp]
    \centering
    \includegraphics[width=0.9\linewidth]{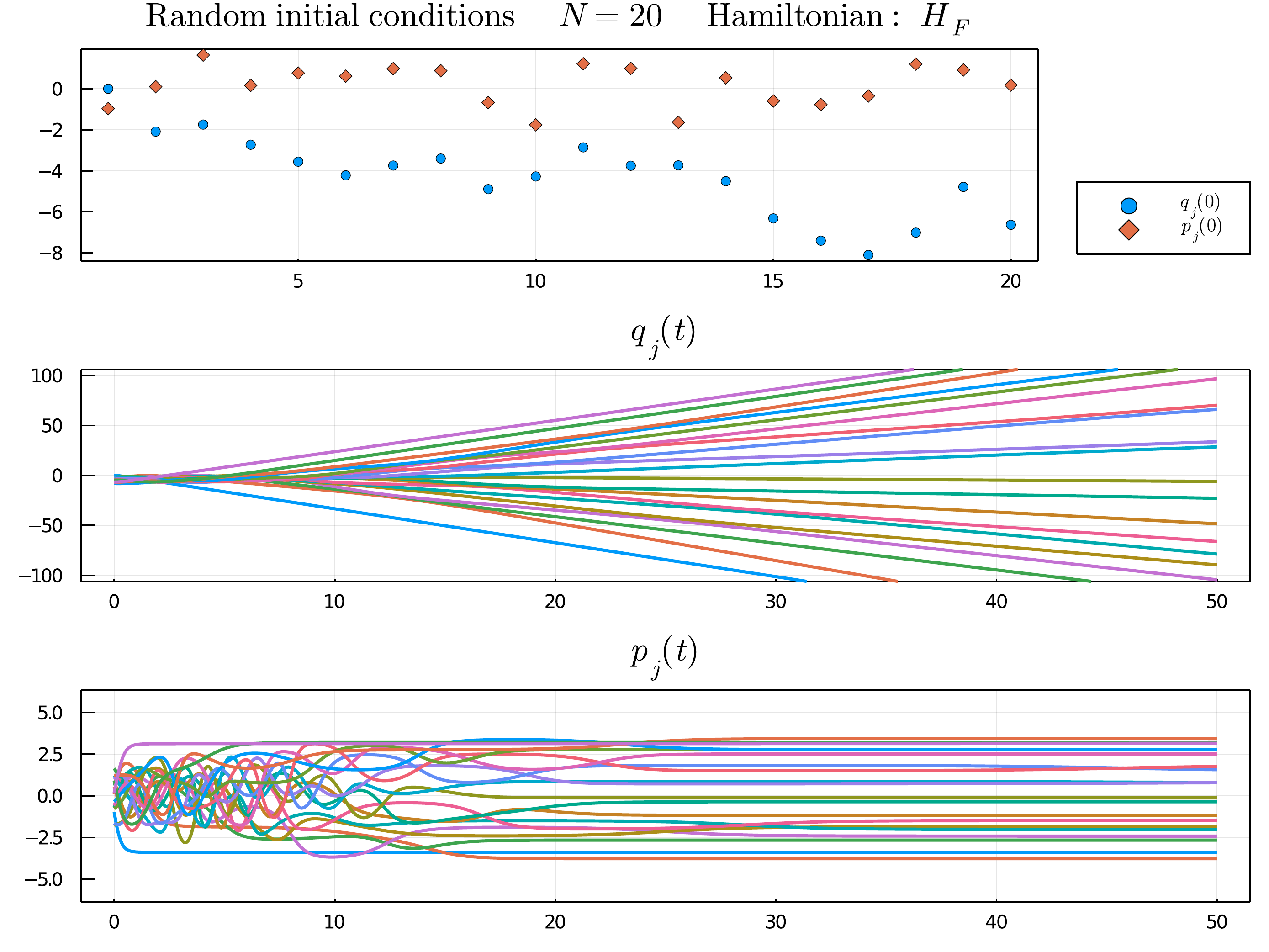}
    \caption{Numerical computations for the open Toda lattice.  The system is integrated using a second-order accurate St\"ormer-Verlet method \cite{BiTr2019} with a time step of $\Delta t = 0.0001$.  The initial conditions are generated by sampling $p_j(0), q_{j+1}(0) - q_j(0)$, $j = 1,\ldots,N$ as independent and normally distributed random variables. }
    \label{f1}
\end{figure}

\begin{figure}[tbp]
    \centering
    \includegraphics[width=0.9\linewidth]{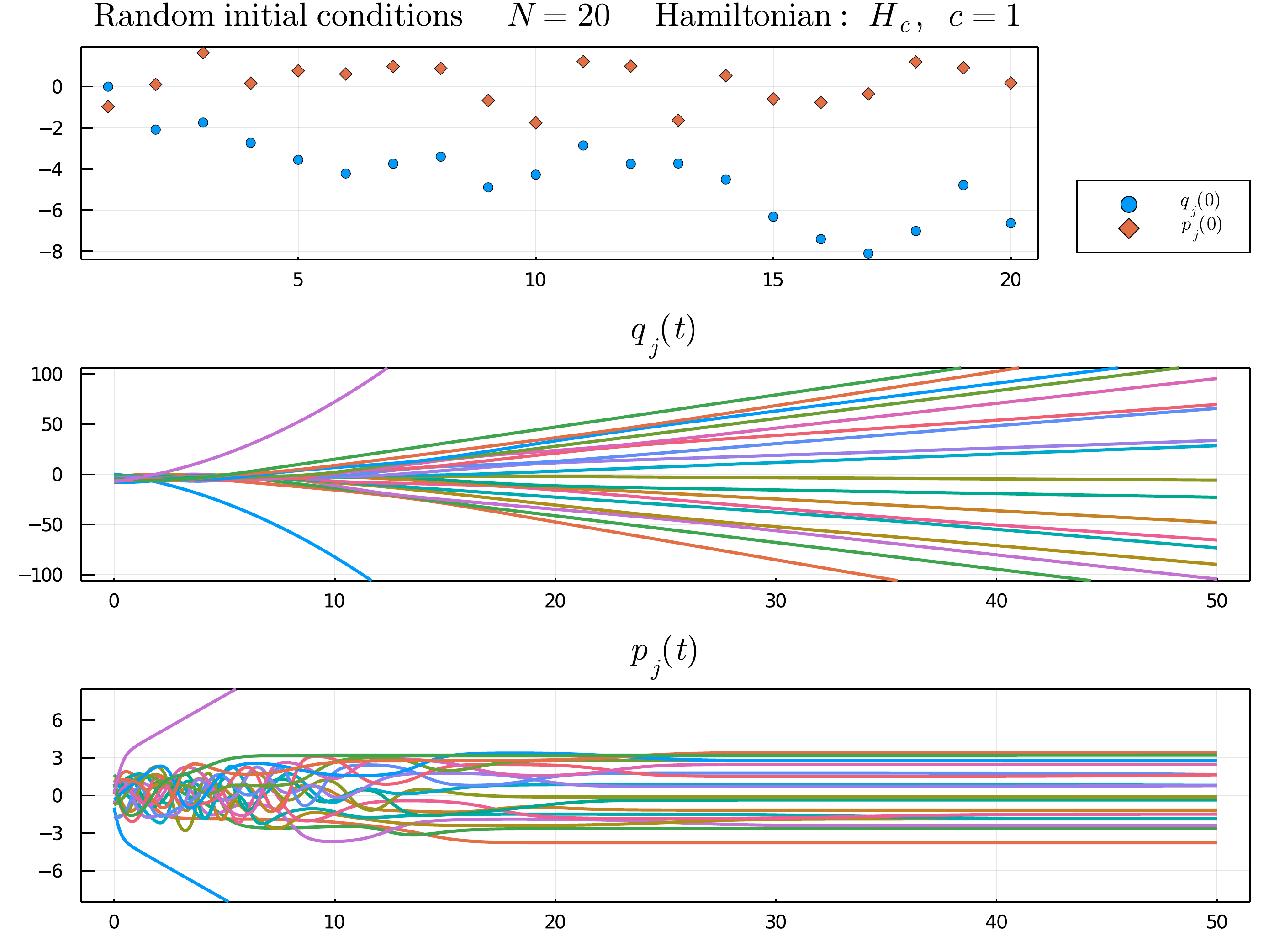}
    \caption{Numerical computations for the perturbed Toda lattice with $c = 1$.  The system is integrated using a second-order accurate St\"ormer-Verlet method with a time step of $\Delta t = 0.0001$.  The initial conditions are generated by sampling $p_j(0), q_{j+1}(0) - q_j(0)$, $j = 1,\ldots,N$ as independent and normally distributed random variables.}
    \label{f2}
\end{figure}

As a benchmark, Figure~\ref{f1} displays the solution of the Toda lattice with Hamiltonian $H_F$, $N=20$ particles and randomly chosen initial data. As $t \to \infty$, $p(t) = p_\infty + o(1)$ and $q(t) = q_\infty + p_\infty t + o(1)$ for suitable constants $q_\infty$ and $p_\infty$ as in (\ref{9.1}) and (\ref{13.1}). Figure~\ref{f2} displays the solution of the perturbed Toda lattice with Hamiltonian $H_c$, $c=1, N=20$ particles and randomly chosen initial data. As $t \to \infty$,
\begin{equation*}
p_i(t) \e p_{i,\infty} + o(1) \ , \quad q_i(t) \e q_{i,\infty} + t p_{i,\infty} + o(1) \ , \quad 2 \le i \le N-1
\end{equation*}
for suitable constants $p_{i,\infty}, q_{i,\infty}$. But
\begin{equation*}
\begin{array}{ll}
p_1(t) \e - c t + O(1) \ , \quad q_1(t) \e - c t^2/2 + O(t) \ , \\
p_N(t) \e c t + O(1) \ , \quad q_N(t) \e c t^2/2 + O(t) \ .
\end{array}
\end{equation*}
This suggests that the solutions of the $H_c$ equations behave like solutions of a system of $N$ particles $q_1, \ldots, q_N, p_1, \ldots, p_N$ 
consisting of a Toda lattice $q_2,\cdots, q_{N-1}$, $p_{2},\cdots, p_{N-1}$ decoupled from a pair of (decoupled) particles $q_1,p_1,$ $q_N, p_N$
solving
\begin{equation*}
\begin{array}{ll}
\dot{p}_1 \e - c \ , \quad \dot{q}_1 \e p_1 \\
\dot{p}_N \e c \ , \quad \dot{q}_N \e p_N \ .
\end{array}
\end{equation*}
Such a system of $N$ particles is clearly completely integrable. What we will show is that solutions of the perturbed Toda system with Hamiltonian $H_c, c>0$, indeed behave asymptotically like solutions of the decoupled system, and hence in view of Moser's observation, the perturbed system is integrable.

\begin{figure}[tbp]
    \centering
    \includegraphics[width=0.9\linewidth]{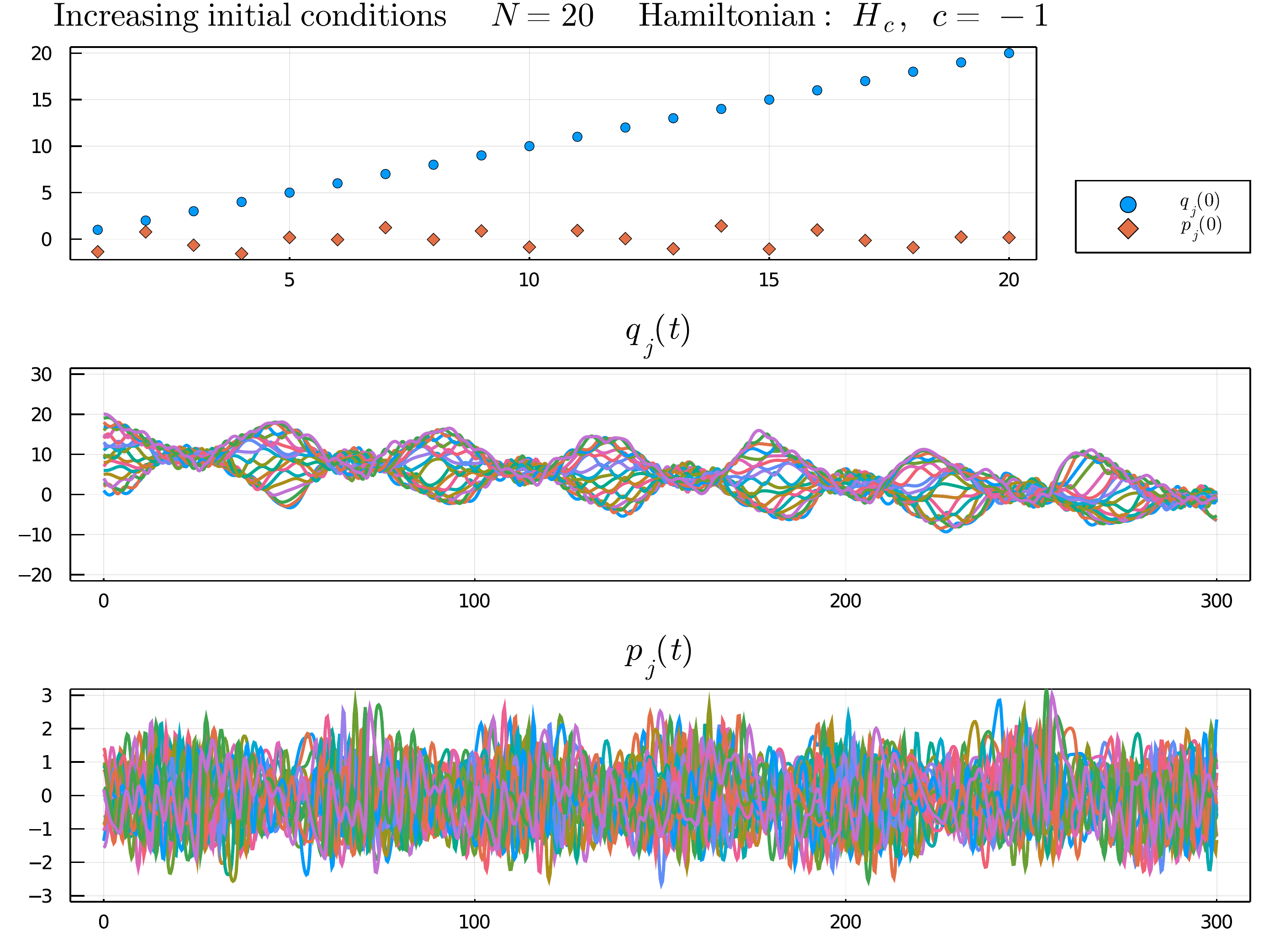}
    \caption{Numerical computations for the perturbed Toda lattice with $c = -1$.  The system is integrated using a second-order accurate St\"ormer-Verlet method with a time step of $\Delta t = 0.0001$.  The initial conditions are generated, for $j =1,\ldots,N$, by sampling $p_j(0)$ as independent and normally distributed and setting $q_{j+1}(0) - q_j(0) = 1$.}
    \label{f3}
\end{figure}

\begin{figure}[tbp]
    \centering
    \includegraphics[width=0.9\linewidth]{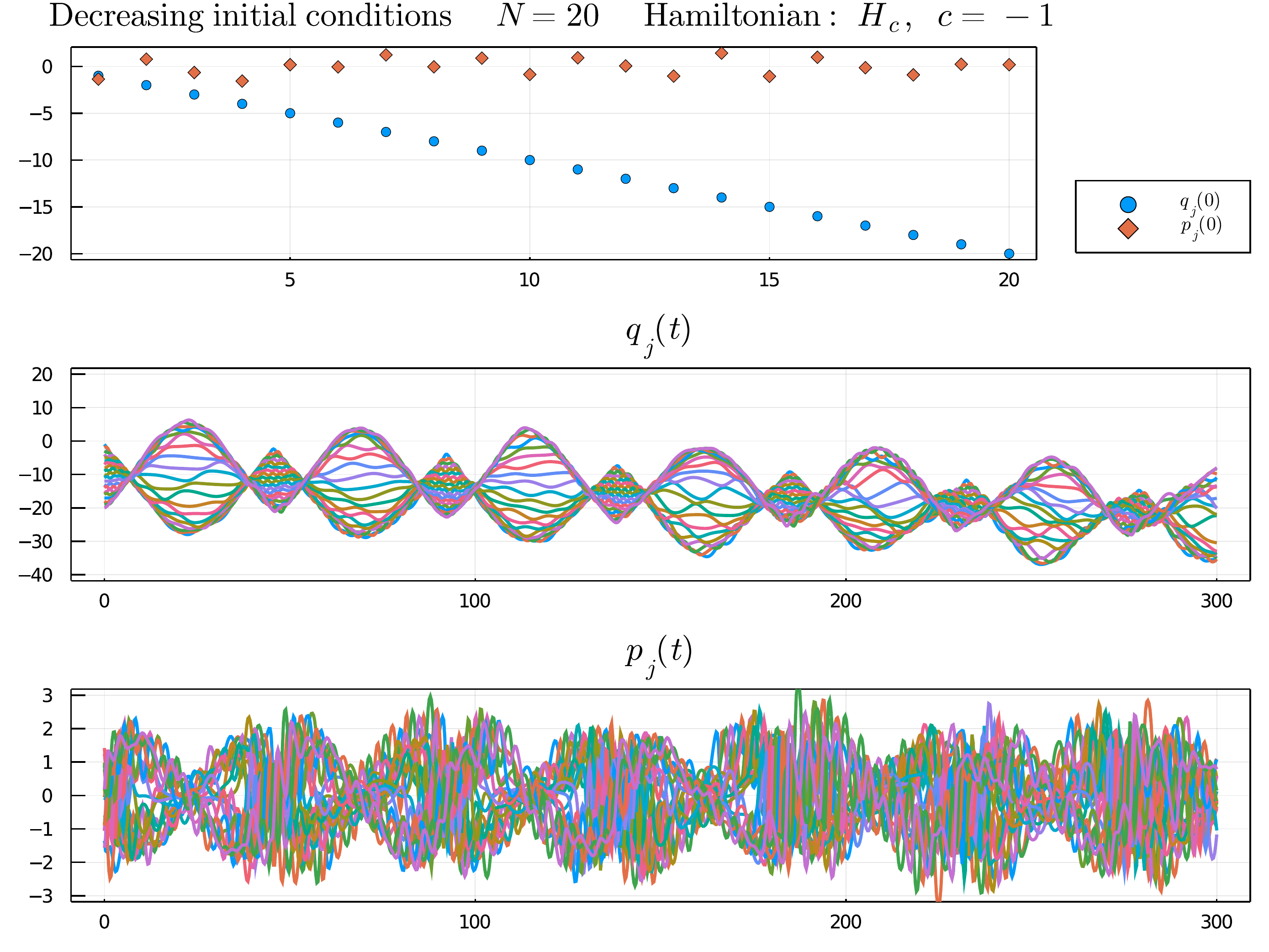}
    \caption{Numerical computations for the perturbed Toda lattice with $c = -1$.  The system is integrated using a second-order accurate St\"ormer-Verlet method with a time step of $\Delta t = 0.0001$.  The initial conditions are generated, for $j =1,\ldots,N$, by sampling $p_j(0)$ as independent and normally distributed and setting $q_{j+1}(0) - q_j(0) = -1$.}
    \label{f4}
\end{figure}

\begin{figure}[tbp]
    \centering
    \includegraphics[width=0.9\linewidth]{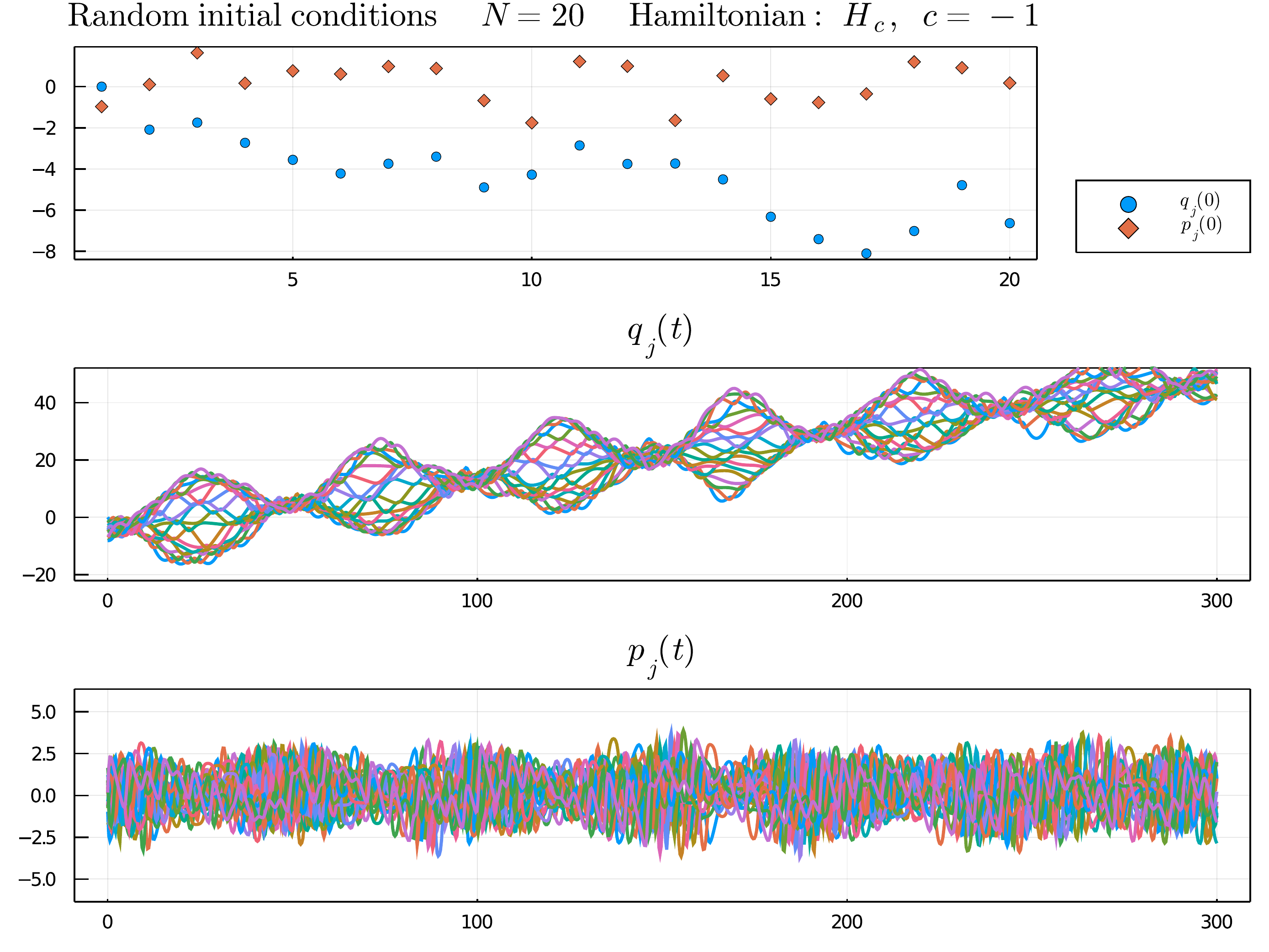}
    \caption{Numerical computations for the perturbed Toda lattice with $c = -1$.  The system is integrated using a second-order accurate St\"ormer-Verlet method with a time step of $\Delta t = 0.0001$.  The initial conditions are generated by sampling $p_j(0), q_{j+1}(0) - q_j(0)$, $j = 1,\ldots,N$ as independent and normally distributed random variables.}
    \label{f5}
\end{figure}

Figure~\ref{f3} displays the solution of the perturbed Toda lattice with Hamiltonian $H_c$ and $c=-1$. Here there are $N$ particles and the initial conditions are
\[ q_i(0) = i \ ,\quad 1 \le i \le 20 \ ,\]
with \{$p_i(0)$\} random. In Figure~\ref{f4}, we again have $c=-1$, but now
\[ q_i(0) = -i \  , \quad 1 \le i \le 20 \ ,\]
again with \{$p_i(0)$\} random.
In Figure~\ref{f5}, we again have $c=-1$, but now \{$q_i(0)$\} and \{$p_i(0)$\}  are chosen randomly. In all three cases, the solution $q_i(t)$ appears to evolve almost periodically in time, modulo a slight gradient. In the
first two cases, this behavior persists at least up to times   $t\approx 300$, but in the third case the almost periodicity begins to unravel after  $t\approx 200$. 

This brings to mind the celebrated computations of E. Fermi, J. Pasta, S. Ulam and M. Tsingou \cite{FePaUlTs1955}, in which the authors, anticipating ergodicity, found, unexpectedly, almost periodic behavior in the solutions of a particular nonlinear lattice system. This meant that in some sense the system was `remembering' its past, and the only way a mechanical system can `remember' its past is if it has many integrals of the motion. In this way, the discovery was viewed as strong evidence for integrability and led eventually, and famously, to the discovery by Kruskal-Zabusky and Gardner-Greene-Kruskal-Miura that the Korteweg de Vries equation is completely integrable.

Over the years, as the power of computers grew, it became clear that Fermi et al. had just not run their equations long enough: With longer computations, they would have found that the almost periodicity unravelled and ergodicity emerged. A very interesting understanding of Fermi et al. is given in \cite{GaPoRi2020}: The lattice equations for unidirectional lattice waves 
can be written schematically in the form
$$ \dot{x} = V(x) +O(h^2) $$

\noindent where $h^2$ is a continuum limit parameter, $h^2 \to 0$, and
$$\dot{y} = V(y)$$

\noindent is KdV. It follows that the solution of the lattice equation $x(t)$ behaves like the (integrable) KdV equation for times $T$ of order  
$h^{-2}$, i.e., $T h^2 $ = O(1), when $x(t)$ begins to diverge from $y(t)$. Thus, the lattice has many $h^2$-accurate integrals
up to times of order $h^{-2}$. It turns out, however,  that the near-integrability
persists for much longer  times $T$ of order $h^{-4}$, and this they are elegantly able to explain by showing that in fact $x(t)$ solves a system of the form
$$ \dot{x} = W(x) +O(h^4) $$
and now 
$$ \dot {y}= W(y)$$
is a solution of the KdV hierarchy, and hence, also, integrable. 

We are led to the following speculation: Is the Fermi et al. problem a guide to what we see for $c<0$? In Figure 5, in particular, when the almost periodicity unravels on a moderate time scale (the same is likely true regarding Figure 3 and Figure 4, but on a longer time scale), Fermi et al. raises the issue 
of whether there is some integrable system associated with the lattice, which describes the solutions of the lattice equations to high accuracy for large, but not infinite, times? In this way, for large times, the system would have excellent, but not perfect, `memory'.

One final comment: The Fermi-Pasta-Ulam-Tsingou paradox, as it is called, is a modern illustration of the interesting phenomenon that sometimes science makes progress, not because of the accuracy of its instruments, but rather because of their inaccuracy. If computers in the 1950's could have made longer calculations, would KdV have been discovered as an integrable system? If Copernicus had more accurate instruments, sensitive to the fluctuations in the planetary orbits, would Kepler have been able to come up with his perfect laws?

\bs
We will prove the integrability of $H_c$ with $c>0$ in steps. In Step 1, we prove that solutions of (\ref{14.2}) with initial data $q_i(0), p_i(0), 1 \le i \le N$ are unique and exist globally, both for $c \ge0$ and $c <0$. In the remainder of the paper, we only consider the case $c>0$. In Step 2, we will show that, as $t \to \infty$, the particle system under $H_c$ splits up into two parts: a core Toda lattice $q_2, \ldots, q_{N-1}, p_2, \ldots, p_{N-1}$ obeying (\ref{7.2}) up to super-exponentially small errors,
\begin{equation*}
\begin{array}{llll}
\dot{q}_n = p_n, \ \quad 2 \le n \le N-1, \\
\dot{p}_2 \e - e^{q_2 - q_3} + O_2(t), \\
\dot{p}_n \e e^{q_{n-1} - q_n} - e^{q_n - q_{n+1}},  \quad 3 \le n \le N-2, \\
\dot{p}_{N-1} \e e^{q_{N-2} - q_{N-1}} - O_{N-1}(t),
\end{array}
\end{equation*}
where $O_2(t) = e^{q_1 - q_2} \e O(e^{-\gamma t^2}), O_{N-1}(t) = e^{q_{N-1} - q_N} \e O(e^{- \gamma t^2})$ for some $\gamma >0$, and a pair of decoupled particles $q_1, q_N, p_1, p_N$ separating from the core lattice, $q_1(t) \to -\infty$ and $q_N(t) \to \infty$, as $t \to \infty$. In Step 3, for solutions $q_1(t), q_2(t), \ldots, q_{N}(t), p_1(t), p_2(t), \ldots, p_N(t)$ of (\ref{14.2}), we obtain precise asymptotics for the inner core $q_2(t), \ldots, q_{N-1}(t), p_2(t), \ldots, p_{N-1}(t)$. 
 Let $U_t(q_0.p_0) = ((q(t),p(t))$ denote the solution of the equations  (\ref{14.2}) generated by $H_c$. Let $\hat U_t(\hat q_0, \hat p_0) = (\hat q(t),\hat p(t))$ denote the solution generated by
\begin{equation*}
H^d_c(q,p) \e \frac{1}{2} \sum_{n=1}^N p_n^2 + \sum_{n=2}^{N-2} \ e^{q_n - q_{n+1}} + c \ \sum_{n=1}^{N-1} (q_n - q_{n+1}),
\end{equation*}
in which the inner Toda core $(q_2,...,q_{N-1},p_2,...,p_{N-1})$ is decoupled from 
particles $q_1$ and $q_N$. Finally, let $U^\#_t (q^\#_0, p^\#_0) = (q^\#(t),p^\#(t))$
denote the solution of the equation generated by the ``free" decoupled Hamiltonian

\begin{equation*}
H^\#_c(q,p) \e \frac{1}{2} \sum_{n=1}^N p_n^2 + c(q_1-q_N).
\end{equation*}
Then in Step 4 we use the asymptotics obtained in Step 3 to show that as $t \to \infty$, solutions of  (\ref{14.2}) behave like ``free'' particles, and the convergence is sufficiently rapid so that Moser's argument applied and the wave operator
$$W^\#(q_0.p_0)  \e \lim_{t \to \infty}  U^\#_{-t} \circ U_t(q_0,p_0)$$
exists. On the other hand, standard Toda asymptotics as in (\ref{9.1}) and (\ref{13.1}), also show that as $t \to \infty$, the solution $(\hat q(t),\hat p(t))$
of the equations generated by $H^{d}_c$, also behave like ``free'' particles, and the convergence is sufficiently rapid so that Moser's argument again applied and the wave operator 
$$\hat W^\#(\hat q_0, \hat p_0)  \e \lim_{t \to \infty} U^\#_{-t} \circ \hat U_t(\hat q_0, \hat p_0)$$
exists. A separate argument then shows that $( \hat W ^\#)^{-1}$ exists and a short calculation then shows that 
$$ W= ( \hat W ^\#)^{-1}  W^\#$$
is an intertwining operator for $\hat U_t$ and $U_t$,
$$\hat U_t \circ W \e W \circ U_t$$
and as $H^d_c$ is integrable, the integrability of $H_c$ follows. The intertwining
relation is not enough, however, to show that as $t \to \infty $, the solutions $U_t$ behave,
as advertised above, 
like solutions $\hat U_t$ of the decoupled system: This is proved using a separate argument. 

Note that we do not construct $W$ directly as a wave operator $$\lim_{t \to \infty}
\hat U_{-t} \circ U_t(q_0,p_0).$$
The technical reason for this is discussed at the end of the section, together with a sketch of the argument that is needed to prove the existence of the limit. We leave the details to the interested reader.

 Finally in Step 5 we display $N$  independent, commuting integrals for the $H_c$ flow and show how the flow can be written in Lax-pair form.

\bs
\noi{\bf Step 1.}

\ms
Standard ODE methods show that (\ref{14.2}) has unique local solutions $(q(t), p(t))$ with $(q(0), p(0)) \e (q_0, p_0)$ for which $H_c(q(t), p(t)) \e h_0 \equiv H_c(q_0, p_0)$. Thus
\begin{equation}\label{20.1}
\frac{1}{2} \sum_{n=1}^N p_n^2(t) + \sum_{n=1}^{N-1} e^{q_n(t) - q_{n+1}(t)} = h_0 - c(q_1(t) - q_N(t))
\end{equation}
and so in order to prove global existence it is enough to show, in particular, that
\begin{equation*} 
| q_1(t) - q_N(t)| \le c_1 t^2 + c_2
\end{equation*}
for some constants, $c_1, c_2$. Indeed, by (\ref{20.1}), we would then have, for $1 \le n \le N$, $| p_n(t)| \le c'_1 t + c'_2$ for some constants $c'_1, c'_2$ and so
$| q_n(t)|\ \le c''_1 t^2 + c''_2$, again for some constants $c''_1, c''_2$. Global existence for $(q(t), p(t))$ now follows, again by standard ODE methods. We derive stronger bounds on $q_n(t), p_n(t)$ in Step 2 below.

\ms
The following elementary calculation plays a crucial role in our analysis. From (\ref{20.1}),
\[ (p_1 - p_N)^2 \le 2 (p_1^2 + p_N^2) \le 2 \sum_{n=1}^{N} p_n^2 \le 4 (h_0 - c(q_1 - q_N)) \ . \]
Setting $q_1 - q_N = \Delta$, we have $(\dot{\Delta})^2 \le 4 (h_0 - c \Delta)$ and so
\begin{equation}\label{21.2}
- 2 \sqrt{ h_0 - c \Delta} \le \dot{\Delta}  \le 2 \sqrt{h_0 - c \Delta} \ .
\end{equation}
Integrating we find
\[ - t + c_3 \le \frac{\sqrt{h_0 - c\Delta}}{-c} \le t + c_4 \]
for some constants $c_3, c_4$. For $c>0$, this implies
\begin{equation}\label{22.0}
0 \le h_0 - c \Delta \le (ct + c')^2
\end{equation}
where $c'$ is some constant and so
\begin{equation*}
\frac{ h_0 - (ct + c')^2}{c} \le q_1 - q_N \le \frac{h_0}{c} \ .
\end{equation*}
There are, of course, similar bounds for $c < 0 $.

\bs
\noi{\bf Step 2.}

\ms
From (\ref{14.2}), we have
\[ \frac{d}{dt} (p_N - p_1) \e e^{q_{N-1} - q_N} + e^{q_1 - q_2} + 2c \]
which implies
\begin{equation}\label{22.2}
p_N(t) - p_1(t) \e p_N(0) - p_1(0) + \int_0^t (e^{q_{N-1} - q_N} + e^{q_1 - q_2}) \ ds + 2 c t .
\end{equation}
Now, from (\ref{21.2}) and  (\ref{22.0}),
\begin{equation} \label{22.3}
p_N - p_1 \e (\dot{q}_N - \dot{q}_1) \le 2 \sqrt{ h_0 + c(q_N - q_1)} \le 2 c t + 2 c' \ .
\end{equation}
Inserting (\ref{22.2}) into (\ref{22.3}), we conclude that
\begin{equation}\label{22.4}
\int_0^\infty (e^{q_{N-1} - q_N} + e^{q_1 - q_2}) \ ds < \infty \ .
\end{equation}
By (\ref{14.2}),
\[ p_1(t) \e p_1(0) - \int_0^t e^{q_1 - q_2} \ ds - c t \]
and so, as $t \to \infty$,
\begin{equation}\label{22.5}
p_1(t) \e p_{1, \infty} - ct + o(1)
\end{equation}
for some constant $p_{1, \infty}$.  Similarly, as $t \to \infty$,
\begin{equation}\label{22.6}
p_N(t) \e p_{N,\infty} + ct + o(1)
\end{equation}
for some constant $p_{N,\infty}$ (note that $p_{i,\infty} \ne \lim_{t \to \infty} p_i(t),$ $i=1,N$). We have
\begin{equation}\label{23.1}
\frac{1}{2} p_1(t)^2 \e \frac{1}{2} c^2\ t^2 - c\ t\ p_{1, \infty} + o(t)
\end{equation}
and
\begin{equation}\label{23.2}
\frac{1}{2} p_N(t)^2 \e \frac{1}{2} c^2\ t^2 +c\ t\ p_{N, \infty} + o(t) \ .
\end{equation}
Inserting (\ref{23.1}) and (\ref{23.2}) into
\[ \frac{1}{2} \sum_{n=1}^N p_n^2 \le h_0 + c\ (q_N - q_1), \]
and using (\ref{22.3}), we find
\[ c^2 \ t^2 + c\ t\ ( p_{N,\infty} - p_{1, \infty}) + \frac{1}{2} \sum_{n=2}^{N-1} p_n^2 + o(t) \]
\[ \le h_0 + c (q_N- q_1) \le c^2\ t^2 + 2 c \ c' t + (c')^2 \ . \]
We conclude that
\beql{23.3}
| p_n| \e O(t^{1/2}) \ , \quad 2 \le n \le N-1
\eeq
and so
\beqn
|q_n| \le O(t^{3/2}) \ , \quad 2 \le n \le N-1 \ .
\eeqn
These bounds are sharper than those obtained in Step 1, and as
\beqn
q_1(t) \e - c \frac{t^2}{2} + O(t) \ , \quad q_N(t) \e c \frac{t^2}{2} + O(t)
\eeqn
by (\ref{22.5}) and (\ref{22.6}), we see that the particles $q_1$ and $q_N$ separate from the core $q_2, \ldots, q_{N-1}$. Moreover,
\beql{24.1}
e^{q_1 - q_2} \ , e^{q_{N-1} - q_N} \e e^{ - c \frac{t^2}{2} (1 + O(t^{-1/2}))} \e O(e^{-\gamma t^2})
\eeq
for some $\gamma >0$.

\bs
\noi{\bf Step 3.}

\ms
Consider $H_F(q_2, \ldots, q_{N-1}, p_2, \ldots, p_{N-1}) \e \frac{1}{2} \sum_{n=2}^{N-1} p_n^2 + \sum_{n=2}^{N-2} e^{q_n - q_{n+1}}$, evaluated  along the solutions $q_1(t), q_2(t), \ldots, q_{N-1}(t), q_N(t),p_1(t), p_2(t), \ldots, p_{N-1}(t), p_N(t)$ of (\ref{14.2}), the flow induced by $H_c$. Then
\[ \frac{d}{dt} H_F(q_2, \ldots, q_{N-1}, p_2, \ldots, p_{N-1}) \e
\sum_{n=2}^{N-1} p_n \dot{p}_n + \sum_{n=2}^{N-2}  e^{q_n - q_{n+1}} (p_n - p_{n+1})\]
\[ \e \sum_{n=2}^{N-1} p_n\  e^{q_{n-1} - q_{n}} - \sum_{n=2}^{N-1} p_n\ e^{q_n - q_{n+1}}\]
\[ + \sum_{n=2}^{N-2} p_n \ e^{q_n - q_{n+1}} -  \sum_{n=2}^{N-2} p_{n+1}\ e^{q_n - q_{n+1}}\]
\[ \e - p_{N-1}\ e^{q_{N-1} - q_{N}} + p_2\ e^{q_1 - q_2} \ . \]

It follows from (\ref{23.3}) and (\ref{24.1}) that
\beql{24.2}
H_F(q_2, \ldots, q_{N-1}, p_2, \ldots, p_{N-1}) \le \hbox{const.}
\eeq
and hence
\beql{24.3}
| p_n(t)| \le \hbox{const.}, \quad 2 \le n \le N-1
\eeq
which is a further strengthening of (\ref{23.3}).

The argument now follows in analogy with Moser's convergence argument for the Toda lattice in \cite{Mo1975}. From (\ref{14.2}),
\[ \frac{d}{dt} p_2 \e e^{q_1 - q_2} - e^{q_2- q_3},\]
and we obtain
\[ p_2(t) \e p_2(0) - \int_0^t e^{q_2 - q_3}\,ds + \int_0^t e^{q_1 - q_2}\,ds \]
\[ \e p_2(0) - \int_0^t e^{q_2 - q_3}\,ds + \int_0^\infty e^{q_1 - q_2} \, ds + O(e^{-\gamma t^2})\]
by (\ref{24.1}). Hence, as $|p_2(t)|$ is bounded by (\ref{24.3}), we conclude that
\beqn
\int_0^\infty e^{q_2 - q_3} \ ds < \infty \ .
\eeqn
Thus
\beql{25.2}
p_2(t) \e p_{2,\infty} + o(1) \ .
\eeq
Now
\[ \frac{d}{dt} (p_2 + p_3) \e e^{q_1 - q_2} - e^{q_3 - q_4}\]
and as $p_2(t)$ and $p_3(t)$ are bounded, we conclude again as above that
\beqn
\int_0^\infty e^{q_3 - q_4} \ ds < \infty
\eeqn
and, using (\ref{25.2}),
\beqn
p_3(t) \e p_{3,\infty} + o(1) \ .
\eeqn
In particular, we conclude that
\beql{26.1}
q_2(t) - q_3(t) \e t (p_{2,\infty} - p_{3,\infty}) + o(t) \ .
\eeq
Using
\[ \frac{d}{dt}(p_2 + \ldots + p_n) \e e^{q_1- q_2} - e^{q_n - q_{n+1}}\ , \quad 2 \le n \le N-1 \]
and proceeding by induction we find
\beql{26.2}
\int_0^\infty e^{q_n - q_{n+1}} \ ds < \infty , \quad 2 \le n \le N-2 \ .
\eeq
(Of course the estimates 
\[ \int_0^\infty e^{q_{N-1} - q_N} \ ds,\, \int_0^\infty e^{q_1 - q_2} \ ds < \infty\]
were obtained earlier.)

It follows that, as $t \to \infty$,
\beqn
p_n(t) = p_{n,\infty} + o(1), \quad 2 \le n \le N-1
\eeqn
and hence
\beqn
q_n(t) \e q_{n,\infty} + p_{n,\infty} t + o(t),  \quad 2 \le n \le N-1
\eeqn
and so
\beql{26.5}
q_n(t) - q_{n+1}(t) = q_{n,\infty} - q_{n-1, \infty} + (p_{n,\infty} - p_{n+1,\infty}) t + o(t), \quad 2 \le n \le N-2 \ .
\eeq
In particular, by (\ref{26.2}), we must have
\beql{26.6}
p_{n,\infty} - p_{n+1,\infty} \le 0, \quad 2 \le n \le N-2 \ .
\eeq
We will  show shortly that the inequality in (\ref{26.6}) is strict. But note first that by (\ref{24.2}) and (\ref{24.3}),
\[ \frac{d}{dt} e^{q_n - q_{n+1}} \e e^{q_n - q_{n+1}} (p_n - p_{n+1}), \quad 2 \le n \le N-2\]
is bounded, and so $e^{q_n - q_{n+1}}$ is globally Lipschitz in time, and in particular uniformly continuous in $t$, and it follows from (\ref{26.2}) that
\beql{27.1}
e^{q_n(t) - q_{n+1}(t)} \to 0
\eeq
pointwise as $t \to \infty$.

\ms
We now prove that the inequality in (\ref{26.6}) is strict.

\ms
In terms of the variables introduced  in (\ref{7.5}),  $a_n \e - p_n/2, 1 \le n \le N$, and $b_n \e e^{(q_n - q_{n+1})/2}/2,\  1 \le n \le N-1$, the equations (\ref{14.2}) take the form
\beql{28.1}
\begin{array}{ll}
\dot{a}_n \e 2 (b_n^2 - b_{n-1}^2), \quad 2 \le n \le N-1, \\
\dot{b}_n \e b_n (a_{n+1} - a_n), \quad 2 \le n \le N-2 \ .
\end{array}
\eeq
Set $A_n = \lim_{t \to \infty} a_n \e - p_{n,\infty}/2 \ , 2 \le n \le N-1$, and consider
\beqn
S(t) \equiv \sum_{n=2}^{N-1} (a_n - A_n)^2 + 2 \sum_{n=2}^{N-2} b_n^2 \ .
\eeqn
From (\ref{28.1}), we have
\[ b_n(t) \e b_n(0) \ e^{\int_0^t (a_{n+1} - a_n) \ ds} \ , \quad 2 \le n \le N-2\]
and as $b_n(t) = \frac{1}{2} e^{(q_n(t) - q_{n+1}(t))/2} \to 0$ by (\ref{27.1}), we must have
\beql{28.2}
\int_0^\infty ( a_{n+1} - a_n) \ ds \e - \infty \ , \quad 2 \le n \le N-2.
\eeq
Using (\ref{28.1}), we find after summing by parts
\[ \frac{1}{4} \frac{d}{dt} S(t) \e \sum_{n=2}^{N-1} (a_n - A_n)(b_n^2 - b_{n-1}^2)
 + \sum_{n=2}^{N-2} b_n^2(a_{n+1}- a_n)\]
\[ \e (a_{N-1} - A_{N-1}) \ \sum_{n=2}^{N-1} (b_n^2 - b_{n-1}^2)\]
\[ + \sum_{n=2}^{N-2} (a_n - a_{n+1} - A_n + A_{n+1}) \ \sum_{i=2}^{n}(b_i^2 - b_{i-1}^2) + \ \sum_{n=2}^{N-2} b_n^2(a_{n+1} - a_n) \]
\[\e (a_{N-1} - A_{N-1})(b_{N-1}^2 - b_1^2)\]\[ + \sum_{n=2}^{N-2} (a_n - a_{n+1} - A_n + A_{n+1})(b_n^2 - b_1^2)+\sum_{n=2}^{N-2} b_n^2(a_{n+1} -a_n)\]
\[=  \sum_{n=2}^{N-2} (A_{n+1} - A_n) b_n^2 + O(e^{-\gamma t^2}),\]
where we have used (\ref{24.1}),
\[ b_1^2 \e \frac{1}{4}\  e^{q_1 - q_2} \ , \quad b_{N-1}^2 \e \frac{1}{4}\ e^{q_{N-1} - q_N} \e  O(e^{-\gamma t^2})\ .\]
Now, if $A_{n+1} - A_n \e - \frac{1}{2}(p_{n+1, \infty} - p_{n,\infty}) > 0$, then it follows from (\ref{26.5}) that
\[ b_n(t) \e \frac{1}{2} e^{(q_n(t) - q_{n+1}(t))/2} \e O(e^{-\mu t /2}) \]
for some $\mu >0$. As $A_{n+1} - A_n \ge0$ by (\ref{26.6}), it follows that
\[ (A_{n+1} - A_n) b_n^2 (t) \e O(e^{-\mu t }), 2 \le n \le N-2 \ . \]
We conclude that, as $ S(t) \to 0 $ as $t \to \infty$,
\[ S(t) \e O(e^{-\mu t })\]
and hence
\[ a_n \e A_n + O(e^{-\mu t }), \quad 2 \le n \le N-1 \]
as $t \to \infty$. Thus if $A_{n+1} - A_n =0$ for some $2 \le n \le N-2$, we would have
$a_{n+1}(t) - a_n(t) = O(e^{-\mu t })$ for some $2 \le n \le N-2$, which contradicts (\ref{28.2}). Thus the inequality in (\ref{26.6}) is strict,
\beql{29.1}
p_{2,\infty} < p_{3, \infty} < \ldots < p_{N-1, \infty} \ .
\eeq

In particular, it follows (use (\ref{26.5}), or more directly $S(t) \e O(e^{-\mu t })$) that, as $t \to \infty$,
\beql{30.1}
e^{q_n(t) - q_{n+1}(t)} \equiv O(e^{-\mu t }),  \quad 2 \le n \le N-2
\eeq
and hence, by the above induction argument we also have
\beql{30.2}
p_n(t) \e p_{n,\infty} + O(e^{-\mu t }) \ , \quad 2 \le n \le N-1,
\eeq
and so
\beql{30.3}
q_n \e q_{n,\infty} + t \ p_{n,\infty} + O(e^{-\mu t }) \ , \quad 2 \le n \le N-1.
\eeq

Equations (\ref{25.2}) and (\ref{27.1}) already imply that the core Toda matrix
\begin{equation}\label{30.4}
L_F \e \begin{pmatrix} a_2 & b_2 & \ldots & 0 \\ b_2 & \ddots & \ddots & 0 \\
\vdots & \ddots & \ddots & b_{N-2} \\0 & 0 & b_{N-2} & a_{N-1} \end{pmatrix}
\end{equation}
converges to a diagonal matrix. By (\ref{30.1}) and (\ref{30.2}),  the convergence is exponential.

\ms
\noi{\bf Remark:} If $L_F$ evolved according to the (exact) Toda flow (\ref{8.3}), then, as $t \to \infty$, $p_n(t) \to p_{n,\infty} = - 2 \lambda_n^F, 2 \le n \le N-1$, where the $\lambda_n^F$'s are the eigenvalues of $L_F(t=0)$. As $L_F(t=0)$ is tridiagonal, the $\lambda_n^F$'s are distinct and so the strict inequality in (\ref{26.6}) is immediate. As $L_F(t)$ solves only a perturbed Toda flow, the strict inequality  requires a more subtle analysis, as above.

\ms
\noi{\bf Step 4.}

\ms
From $p_1(t) \e p_1(0) - \int_0^t e^{q_1 - q_2} \ ds - ct$ and (\ref{24.1}), we see that
\beql{31.1}
p_1(t) \e p_{1,\infty} - c t + O(e^{-\gamma t^2})
\eeq
and similarly
\beqn
p_N(t) \e p_{N,\infty} + c t + O(e^{-\gamma t^2})
\eeqn
from which it follows that
\beqn
q_1(t) \e q_{1,\infty} - \frac{c t^2}{2} + t p_{1,\infty} + O(e^{-\gamma t^2})
\eeqn
and
\beqn
q_N(t) \e q_{N,\infty} + \frac{c t^2}{2} + t p_{N,\infty} + O(e^{-\gamma t^2}) \ .
\eeqn
From (\ref{30.2}) and (\ref{30.3}), we have for $2 \le n \le N-1$,
\beqn
p_n(t) \e p_{n,\infty} + O(e^{-\mu t})
\eeqn
and
\beql{31.6}
q_n(t) \e q_{n,\infty} + t p_{n,\infty} + O(e^{-\mu t})\ .
\eeq
Let $U_t(q_0, p_0) \e (q(t), p(t))$ be the solution of (\ref{14.2}) with $(q(0),p(0) \e (q_0,p_0)$, and let $\hat U_t(\hat q_0, \hat p_0) \e (\hat q(t), \hat p(t))$ be the solution of the Hamiltonian equations generated by the decoupled Hamiltonian
\beql{31.7}
H_c^d(q,p) \e \frac{1}{2} \sum_{n=1}^N p_n^2 + \sum_{n=2}^{N-2} e^{q_n - q_{n+1}} + c (q_1 - q_N) \ .
\eeq
The Hamiltonian $H_c^d$ is clearly completely integrable. We have
\[\dot{ \hat{q}}_1(t) \e \hat p_1 (t) \ , \quad \dot{\hat{p}}_1 \e - c \]
and so
\beqn
\hat p_1 (t) \e  \hat p_{1,0} - c t \ , \quad \hat q_1(t) \e \hat q_{1, 0} + \hat p_{1,0} t - \frac{ct^2}{2}
\eeqn
and similarly
\beqn
\hat p_N (t) \e \hat p_{N,0} + c t \ , \quad \hat q_N(t) \e \hat q_{N, 0} + \hat p_{N,0} t + \frac{ct^2}{2} \ .
\eeqn
By standard Toda asymptotics, (\ref{9.1}) and (\ref{13.1}), for $2 \le n \le N-1$, as $t \to  \pm \infty$,
\beqn
\hat p_n(t) \e \alpha_n^\pm + O(e^{-\gamma|t|})
\eeqn
and
\beqn
\hat q_n(t) \e \beta_n^\pm + \alpha_n^\pm t +  O(e^{-\gamma|t|}) \ .
\eeqn
Also  let $U_t^\#(q_0^\#,p_0^\#) = (q^\#(t), p^\#(t))$ be the solution of the equations generated by the ``free" decoupled Hamiltonian
\beqn
H_c^\# (q,p) \e \frac{1}{2} \sum_{n=1}^N p_n^2 + c (q_1 - q_n) \ .
\eeqn
Clearly
\beqn
p_1^\#(t) \e p_{1,0}^\# - c t \ , \quad q_1^\#(t) \e q_{1,0}^\# + p_{1,0} ^\# t - \frac{c t^2}{2}
\eeqn
and
\beqn
p_N^\#(t) \e p_{N,0}^\# + c t \ , \quad q_N^\#(t) \e q_{N,0}^\# + p_{N,0} ^\# t + \frac{c t^2}{2}
\eeqn
and
\beql{32.8}
p_n^\#(t) \e p_{n,0}^\#  \ , \quad  q_n^\#(t) \e q_{n,0}^\# + p_{n,0} ^\# t  \ , \quad 2 \le n \le N-1 \ .
\eeq
Now
\[ \hat W_t^\# (\hat q_0, \hat p_0) \equiv U_{-t}^\# \circ \hat U_t (\hat q_0, \hat p_0)\e U_{-t}^\# \big( \hat q_{1,0} + \hat p_{1,0} t - \frac{c t^2}{2}, \hat p_{1,0} - c t,\]
\[ \big( \beta_n^+ +\alpha_n^+ t + O(e^{-\gamma t}), \alpha_n^+ +  O(e^{-\gamma t})\big)_{n=2}^{N-1},  \hat q_{N,0} + \hat p_{N,0} t + \frac{c t^2}{2} , \hat p_{N,0} + c t \big)\]

\[ \e \big( (\hat q_{1,0} + \hat p_{1,0} t - \frac{c t^2}{2}) + (\hat p_{1,0} - c t)(-t) - \frac{c(-t)^2}{2},
 (\hat p_{1,0} - c t) - c(-t)\]

\[  \big( \beta_n^+ + \alpha_n^+ t + O(e^{-\gamma t}) + (\alpha_n^+ +  O(e^{-\gamma t})) (-t), \alpha_n^+ +  O(e^{-\gamma t})\big)_{n=2}^{N-1},
\]
\[(\hat q_{N,0} + \hat p_{N,0} t + \frac{c t^2}{2}) + (\hat p_{N,0} + c t )(-t)+ \frac{c(-t)^2}{2}, (\hat p_{N,0} + c t) + c(-t)\big) \]
\[ \e \big(\hat q_{1,0}, \hat p_{1,0},\big(\beta_N^+ + O(e^{-\gamma t}), \alpha_n^+  + O(e^{-\gamma t})\big)_{n=2}^{N-1}, \hat q_{N,0}, \hat p_{N,0}\big) \]
\[\to  \big(\hat q_{1,0}, \hat p_{1,0}, \big(\beta_N^+ , \alpha_n^+  \big)_{n=2}^{N-1}, \hat q_{N,0}, \hat p_{N,0}\big) \equiv (\hat q^\#, \hat p^\#),\]
where
\beql{33.1}
\hat q^\# \e (\hat q_{1,0}, \beta_2^+, \ldots, \beta_{N-1}^+, \hat q_{N,0}) \ , \quad
\hat p^\# \e (\hat p_{1,0}, \alpha_2^+, \ldots, \alpha_{N-1}^+, \hat p_{N,0})\ .
\eeq
Thus
\beql{33.2} \hat W^\#(\hat q_0, \hat p_0) \e \lim_{t \to \infty} U_{-t}^\# \circ \hat U_t(\hat q_0, \hat p_0) \e (\hat q^\#, \hat p^\#)
\eeq
exists. A similar argument shows that
\beql{33.3}
W^\#(q_0,p_0) \e \lim_{t \to \infty} U_{-t}^\# \circ U_t (q_0,p_0) \e (q^\#, p^\#)
\eeq
exists, where we now use (\ref{31.1}) and (\ref{31.6}),
\beqn
q^\# \e (q_{1,\infty},\ldots, q_{N,\infty}) \ , \quad p^\# \e (p_{1,\infty},\ldots, p_{N,\infty}) \ .
\eeqn
As
\beqn
I_1(q,p) \e \frac{1}{2} p_1^2 + c q_1 \ , \quad I_n(q,p) \e p_n, \ 2 \le n \le N-1 \ , \quad I_N(q,p) \e \frac{1}{2}p_N^2- q_N
\eeqn
are commuting integrals for $U_t^\#$, it follows from (\ref{33.2}), as in (\ref{12.1}), that
\beqn
J_n(q,p) \e I_n \circ W^\#(q,p) \e I_n(q^\#, p^\#) \ , \quad 2 \le n \le N-1
\eeqn
are commuting integrals for $U_t$, i.e., for (\ref{14.2}). In particular, this shows that $H_c$ is integrable. However, we want to show more: we want to show that there is an intertwining operator $W$ for $\hat U_t$ and $U_t$,
\beql{34.4}
\hat U_t \circ W = W \circ U_t \ .
\eeq
This will then allow us to display (\ref{14.2}) in a convenient Lax-pair form, and hence as an isospectral deformation similar to the (unperturbed) Toda case. From (\ref{33.2}) and (\ref{33.3}), we have
\beql{34.5}
U_t^\# \circ \hat W^\# \e \hat W^\# \circ \hat U_t \ \ \hbox{and} \ \
U_t^\# \circ  W^\# \e  W^\# \circ U_t
\eeq
and so, if $(\hat W^\#)^{-1}$ exists, then (modulo domain issues, explained at the end of Step 4)
\[ \hat U_t \circ (\hat W^\#)^{-1} \circ W^\# \e (\hat W^\#)^{-1} \circ U_t^\# \circ W^\# \e (\hat W^\#)^{-1} \circ W^\# \circ U_t\]
and so (\ref{34.4}) holds with
\beql{34.6}
W \e (\hat W^\#)^{-1} \circ W^\# \ .
\eeq
The proof that $(\hat W^\#)^{-1}$ exists requires detailed knowledge of $(\hat q^\#, \hat p^\#)$.

\ms
 To prove that $(\hat W^\#)^{-1}$ exists, we use the following result from \cite{DeDuTr2021} (for an alternative proof, see Props. 5.2 and 3.6 in \cite{LeSaTo2008}). Let
\beql{35.0}
(q(t), p(t)) \e (q_2(t), \ldots, q_{N-1}(t), p_2(t), \ldots, p_{N-1}(t))
\eeq
solve the Toda equation generated by $H_F(q,p) \e \frac{1}{2} \sum_{n=2}^{N-1} p_n^2 + \sum_{n=2}^{N-2} e^{q_n - q_{n+1}}$
with initial data $(q_0, p_0)$. Then, as $t \to \infty$ (see (\ref{9.1})),
\beql{35.1}
q_n(t) \e \alpha_n^+ t + \beta_n^+ + O(e^{-\gamma t}) \ , \quad \gamma >0,\,\, 2 \le n \le N-1,
\eeq
where
\beql{35.2}
\alpha_n^+ \e - 2 \lambda_n \ , \quad 2 \le n \le N-1
\eeq
and
\beql{35.3}
\beta_n^+ \e \frac{1}{N-2} \sum_{j=2}^{N-1} q_{j,0} - \frac{2}{N-2} \sum_{j=2}^{N-1} \ln \left( \frac{u_n(2)}{u_j(2)} \frac{ \Pi_{\ell=2}^{n-1} 2 (\lambda_\ell - \lambda_n)}{\Pi_{\ell=2}^{j-1} 2 (\lambda_\ell - \lambda_j)} \right) \ ,
\eeq
where $2 \le n \le N-1$, and $\Pi_{\ell=2}^1 2 (\lambda_\ell - \lambda_2) \equiv 1$.
Here $\lambda_2 > \lambda_3 > \ldots > \lambda_{N-1}$ are the eigenvalues of the core Toda matrix $L_F(q_0, p_0)$ in (\ref{30.4}), where $a_i \e - p_{i,0}/2$, $b_j \e \frac{1}{2} e^{(q_{j,0} - q_{j+1,0})/2}$, $2 \le i \le N-1$, $2 \le j \le N-2,$ and $u_2(2), \ldots u_{N-1}(2)$ are the first components of the normalized eigenvectors $u_n = (u_n(2), \ldots, u_n(N-1))^T$ of $L_F(q_0,p_0)$ corresponding to $\lambda_n,\ 2 \le n \le N-1$. We have $\sum_{j=2}^{N-1} u_n^2(j) =1$ and $u_n(2)>0$ for all $2 \le n \le N-1$. It is well known (see e.g. \cite{DeDuTr2021}) that the map $\Phi$ from Jacobi matrices $L_F$ with $b_i>0$, $2 \le i \le N-1$, is a bijection onto
\[ \{ (\gamma_2, \ldots,\gamma_{N-1},\mu_2, \ldots, \mu_{N-1}) \ : \ \gamma_2 > \gamma_3 > \ldots > \gamma_{N-1}, \]\[ \sum_{i=2}^{N-1} \mu_i^2=1 , \quad \mu_i >0 , \quad 2 \le i \le N-1 \}. \]
Let $U_t^F \e (q(t), p(t))$ denote the solution of the Toda equations in  (\ref{35.0}) above and let $U_t^0$ denote the solution of the equations generated by $H^0(q,p) \e \frac{1}{2} \sum_{n=2}^{N-1} p_n^2$. Then the argument leading to (\ref{33.2}) shows that
\[ W^F(q_0,p_0) \e \lim_{t \to \infty} U_{-t}^0 \circ U_t^F(q_0, p_0) \]
exists and
\beql{36.1}
W^F(q_0,p_0) \e (\beta_2^+, \ldots, \beta_{N-1}^+, \alpha_2^+,\ldots, \alpha_{N-1}^+) \ .
\eeq
We show first that $W^F$ is one-to-one. Assume (\ref{36.1}). From (\ref{35.2}), the eigenvalues $\{\lambda_n\}$ of $L_F(q_0,p_0)$ are determined,
\beqn
\lambda_n \e - \alpha_n^+/2 \ , \quad 2 \le n \le N-1 \ .
\eeqn
From (\ref{35.3}),
\beql{37.1}
\sum_{n=2}^{N-1} \beta_n^+ \e \sum_{j=2}^{N-1} q_{j,0}  - \frac{2}{N-2} \sum_{n=2}^{N-1} \sum_{j=2}^{N-1} \ln \left( \frac{u_n(2)}{u_j(2)} \frac{ \Pi_{\ell=2}^{n-1} 2 (\lambda_\ell - \lambda_n)}{\Pi_{\ell=2}^{j-1} 2 (\lambda_\ell - \lambda_j)} \right) \e \sum_{i=2}^{N-1} q_{j,0} \ ,
\eeq
as the double sum vanishes by oddness. It follows then from (\ref{35.3}) that, for $2 \le n \le N-1$,
\beql{37.2}
\sum_{j=2}^{N-1} \ln \frac{u_n(2)}{u_j(2)} \e r_n \ ,
\eeq
where $r_n$ is a function of $\{\beta_j^+\}_{i=2}^{N-1}$ and $\{\alpha_n^+ \e - 2 \lambda_n\}$,
\beql{37.3}
r_n \e \frac{N-2}{2}\left( \frac{1}{N-2} \sum_{j=2}^{N-1} \beta_j^+ - \beta_n^+ - \frac{2}{N-2} \sum_{j=2}^{N-1} \ln \left( \frac{ \Pi_{\ell=2}^{n-1} (\alpha_n^+ - \alpha_\ell^+)}{\Pi_{\ell=2}^{n-1} (\alpha_j^+ - \alpha_\ell^+)}\right) \right).
\eeq
Note that
\beql{37.4}
\sum_{n=2}^{N-1} r_n \e 0 \ .
\eeq
From (\ref{37.2}),
\[ \prod_{j=2}^{N-1} \left( \frac{u_n(2)}{u_j(2)}\right) \e e^{r_n}\]
and so
\[ u_n(2) \e \left(\prod_{j=2}^{N-1} u_j(2) \right)^{1/{(N-2)}} \ e^{ r_n/(N-2)} \ . \]
Since $\sum_{n=2}^{N-1} u_n(2)^2 \e 1$, this gives
\[ \prod_{j=2}^{N-1} u_j(2) \e \left( \sum_{n=2}^{N-1} \ e^{2 r_n/(N-2)}\right)^{(2-N)/2} \ . \]
Thus
\beqn
u_n(2) \e \frac{e^{\frac{r_n}{N-2}}}{\big(\sum_{n=2}^{N-1} e^{\frac{2 r_n}{N-2}}\big)^{1/2}} \ , \quad 2 \le n \le N-1 \ .
\eeqn
As $L_F(q_0, p_0)$ is determined by its eigenvalues and the first components of its normalized eigenvectors, it follows that $L_F(q_0, p_0)$ is determined by $\{ \beta_i^+, \alpha_j^+\}, 2 \le i, j, N-1$. But
\[ p_{n,0} \e - 2 a_n, \ \quad 2 \le n \le N-1 \]
and
\[ q_{n,0} - q_{n+1,0} \e 2 \ln 2 b_n \ , \quad 2 \le n \le N-2 \ . \]
As $\sum_{n=2}^{N-1} q_{n,0} \e \sum_{n=2}^{N-1} \beta_n^+$, by (\ref{37.1}),  we thus see that $(q_0,p_0)$ is determined by
\[ (\beta_2^+, \ldots, \beta_{N-1}^+, \alpha_2^+, \ldots \alpha_{N-1}^+) \e W^F(q_0,p_0) \ , \]
i.e., $W^F$ is one-to-one. We now show that $W^F$ is onto
\beqn
X \e \{ (x_2, \ldots, x_{N-1}, y_2, \ldots, y_{N-1}) \ : \ y_2 < \ldots < y_{N-1} \} \subset \RR^{2(N-2)} \ .
\eeqn
Let $r_n = r_n(x,y)$ in (\ref{37.3}) with $\beta_n^+$ replaced by $x_n$ and $\alpha_n^+$ by $y_n,$ $2\leq n\leq N-1.$
Set
\beql{39.1}
u_n(2) \e \frac{e^{r_n/(N-2) }} { \big(\sum_{j=2}^{N-1} \ e^{2 r_j/(N-2)}  \big)^{1/2}} > 0 \ , \quad 2 \le n \le N-1 \ .
\eeq
Then by (\ref{37.4}),
\beqn
\prod_{n=2}^{N-1} u_n(2) \e
\frac{e^{\frac{1}{N-2} \sum_{n=2}^{N-1} r_n}}
 {\big(\sum_{j=2}^{N-1} \ e^{2 r_j/(N-2)}  \big)^{(N-2)/2}} \e \frac{1} {\big(\sum_{j=2}^{N-1} \ e^{2 r_j/(N-2)}  \big)^{(N-2)/2}}\ .
\eeqn
From (\ref{39.1}),
\beql{39.3} \sum_{n=2}^{N-1} u_n(2)^2 \e 1 \quad \hbox{and} \quad u_n(2) \ge 0, \,\,2 \le n \le N-1,
\eeq
and using (\ref{39.1}), we see that
\[ u_n(2) \e \big( \prod_{j=2}^{N-1} u_j(2)\big)^{1/(N-2)} \ e^{r_n/(N-2)}  \]
which implies
\beqn
\prod_{i=2}^{N-1} \ \left( \frac{u_n(2)}{u_j(2)}\right) \e e^{r_n(x,y)}
\eeqn
from which we conclude that
\beqn
x_n \e \frac{1}{N-2} \sum_{j=2}^{N-1} x_j - \frac{2}{N-2} \sum_{j=2}^{N-1} \ln \left(  \frac{u_n(2)}{u_j(2)} \frac{\Pi_{\ell=2}^{n-1} (y_n - y_\ell)} {\Pi_{\ell=2}^{n-1} (y_j- y_\ell)} \right) \ .
\eeqn
Now, as $\Phi$ is a bijection, there exists a unique Toda matrix $L_F$ (see (\ref{30.4})) with spectrum $\lambda_2 = -y_2/2 > \lambda_3 = - y_3/2 > \ldots > \lambda_{N-1} = - y_{N-1}/2$ and first components of the eigenvectors $u_2(2), \ldots, u_{N-1}(2)$.

Set
\beqn
\begin{array}{ll}
p_{n,0} \e - 2 a_n \ , \quad 2\le n \le N-1, \\
q_{n,0} - q_{n+1,0} \e 2 \ln (2 b_n) \ , \quad 2 \le n \le N-1 \ .
\end{array}
\eeqn
Then determine the $q_{n,0}$'s uniquely by requiring
\beqn
\sum_{n=2}^{N-1} q_{n,0} \e \sum_{n=2}^{N-1} x_n \ .
\eeqn
It then follows from the above calculations that
\[ W^F(q_0,p_0) \e (x_2, \ldots, x_{N-1}, y_2, \ldots, y_{N-1} ) \]
which completes the proof that $W^F$ is a bijection from $\RR^{2(N - 2)}$ to $X$. Finally, we conclude from (\ref{33.1}) that $\hat W^\#$ is a bijection from $\RR^{2N}$ to
\[ \hat X^\# \e \{ (x_1, x_2, \ldots, x_{N-1}, x_N, y_1, y_2, \ldots, y_{N-1}, y_N) \ : \ y_2 < y_3 < \ldots < y_{N-1} \} \ .\]
In order to derive (\ref{34.4}) with $W$ as in (\ref{34.6}), we need to verify certain domain issues. For $(x,y) \in \hat X^\#$, we have from (\ref{34.5})
\[ \hat W^\# \circ \hat U_t \circ (\hat W^\#)^{-1} (x,y) \e U_t^\# \circ \hat W^\# \circ (\hat W^\#)^{-1} (x,y) \e U_t^\#(x,y) \]
from which we see necessarily that $U_t^\#(x,y) \in \hat X^\#$, a fact that can be seen directly from (\ref{32.8}). Hence
\[ \hat U_t \circ (\hat W^\#)^{-1} \e (\hat W^\#)^{-1}  \circ U_t^\# \ \ \hbox{on} \ \ \hat X^\# \ . \]
But it follows from (\ref{29.1}) that for any $(x,y) \in \RR^{2N}$,
$W^\#(x,y) \in \hat X^\#$, and so
\[ \hat U_t \circ (\hat W^\#)^{-1} \circ W^\#(x,y) \e (\hat W^\#)^{-1} \circ U_t^\# \circ W^\# \e (\hat W^\#)^{-1} \circ W^\# \circ U_t(x,y) \]
which verifies, indeed, that $W \e (\hat W^\#)^{-1} \circ W^\#$ mapping $\RR^{2N}$ to itself intertwines $U_t$, the propagator for the equations generated by $H_c$, and $\hat U_t$, the propagator for the equations generated by the completely integrable Toda-core Hamiltonian $H_c^d$.

As noted earlier, although the intertwining relation $\hat U_t \circ W = W \circ U_t$ is enough to prove integrability, it is not sufficient to prove that solutions generated by $H_c$ behave asymptotically like solutions generated by the decoupled Toda core Hamiltonian $H^d_c$.  In the quantum mechanical case , the fact that $e^{iAt}$ is linear and  unitary implies that $|| e^{-iAt}e^{iBt} f - W f|| = ||e^{iBt}f - e^{iAt}W f|$, and so the convergence of the wave operator $W$ is equivalent  to showing 
that  a solution $e^{iBt}f$ generated by the operator B, behaves like a solution
$e^{iAt}g$ generated by the operator A, where $g := Wf$. In the case at hand, as $\hat U^\#_t$, is neither linear nor bounded, we cannot, in particular, immediately infer 
from the convergence $U^\#_{-t} \circ \hat U_t(q_0,p_0) - \hat  W^\#(q_0,p_0) \rightarrow 0$, that $\hat U_t(q_0,p_0) - U^\#_t \hat W^\#(q_0,p_0) \rightarrow 0$, as desired. However, as we see from the calculations following (\ref{32.8}), as $t \to \infty$,
$$ U^\#_{-t} \circ \hat U_t(q_0,p_0) = \hat W^\#(q_0,p_0) + O(e^{-\gamma t})$$
it then follows from the explicit form, and polynomial growth,  of $U^\#_t$ that

$$ \hat U_t (q_0,p_0) = U^\#_t (\hat W^\#(q_0,p_0) + O(e^{-\gamma t})) = ( U^\#_t \hat W^\#(q_0,p_0)) + O(e^{-\gamma t/2}).$$  

A similar argument shows that
$$  U_t (q_0,p_0) = U^\#_t ( W^\#(q_0,p_0) + O(e^{-\gamma t})) = ( U^\#_t  W^\#(q_0,p_0)) + O(e^{-\gamma t/2}).$$ 
Now as $W^\#(q_0,p_0) \in \hat X^\#$, and as $\hat W^\#$ is a bijection onto $\hat X^\#$, it follows that there exist $(\hat q_0, \hat p_0)$ such that  $\hat W^\#(\hat q_0,\hat p_0) = W^\#(q_0,p_0).$ Substitution into the above two relations shows that
$$ U_t(q_0,p_0) = \hat U_t(\hat q_0, \hat p_0) + O(e^{-\gamma t/2})$$ 
where $(\hat q_0, \hat p_0) = W (q_0, p_0)$
as desired.

Finally, as noted before, we do not construct $W$ directly as a wave operator $ \lim_{t \to \infty}
\hat U_{-t} \circ U_t(q_0,p_0).$ The reason for this is the following. In evaluating
$$\hat U_{-t} \circ U_t(q_0,p_0) = \hat U_{-t}(q(t),p(t))$$
we are facing a double scaling limit. The asymptotics of $\hat U_{-t}(\hat q_0,\hat p_0)$  as $t \to \infty$ is known for $(\hat q_0, \hat p_0)$ fixed,  or in a compact set, but $ q(t)$, in particular, grows linearly. This considerably complicates the analysis. The difficulty is avoided when we evaluate  $$\hat U^\#_{-t} \circ U_t(q_0,p_0) = \hat U^\#_{-t}(q(t),p(t))$$
as we have an explicit, and simple, formula for $ \hat U^\#_{-t} (\hat q^\#, \hat p^\#)$ for {\bf all} $(\hat q^\#, \hat p^\#)$, and so the double scaling limit is avoided. To avoid the problem of the double scaling limit in evaluating $\hat U_{-t}(q(t),p(t))$, we need to use an explicit formula for
$ \hat U_{-t}( \hat q_0, \hat p_0)$ for all $ (\hat q_0, \hat p_0)$.  This is most
conveniently done by mapping the solution $(q_1(t),...,q_N(t), p_1(t),...,p_N(t))$
onto the eigenvalues $(\lambda_2(t),...,\lambda_{N-1}(t))$ and first components of the associated normalized eigenvectors $(u_2(2)(t), ...,u_{N-1}(2)(t))$ of the core Toda matrix $L_F =L_F(t)$ in (\ref{30.4}) . In \cite{Mo1975}, Moser used the Lax-Pair form (\ref{8.3}) to show that under the Toda flow
$$\hat \lambda_i(t) = \hat  \lambda_{i,0} ,\ i= 2,...,N-1$$ and
$$\hat  u_i(2)(t) =\frac {\hat u_{i,0}(2) e^{\hat \lambda_{i,0}t} }{ (\sum_{j=2}^{N-1} (\hat u_{j,0}(2))^2 e^{2 \hat \lambda_{j,0}t)})^{1/2}},\ i=2,...,N-1$$
for any initial conditions $\hat \lambda_{i,0}$,$\hat u_{j,0}(2)$. The evolution of $L_F(t)$ is almost the same as in (\ref{8.3}), except there is an additional diagonal driving term
$$\dot{L}_F \e [ L_F, B_F] + \diag ( -2b_1^2,0,0,...,0,2b_{N-1}^2)$$
where $b_1= \frac{1}{2}e^{(q_1(t)-q_2(t))/2}$ and $b_{N-1}= \frac{1}{2}e^{(q_{N-1}(t)-q_N(t))/2}.$ Now Moser's method to obtain the above formulae for ($\hat \lambda_i(t)$,$\hat u_j(2)(t)$), can be extended, using the driven Lax-Pair
for $L_F(t)$, to obtain the asymptotics of $\lambda_i (t)$ and $u_j(2)(t)$. Then when we evaluate
$$\hat U_{-t} \circ U_t(q_0,p_0) = \hat U_{-t}(q(t),p(t))$$
now in the ($\lambda, u(2)$) variables, rather than the original $(q,p)$ variables, we can use Moser's explicit formulae for ($\hat \lambda_i(t)$,$\hat u_j(2)(t)$), and the double scaling limit is avoided. Using the bijection between $L_F$  and the ($\lambda, u(2)$)
variables, we can then assert, after some algebra,  the existence of
$$ \lim_{t \to\infty}\hat U_{-t} \circ U_t(q_0,p_0)$$ directly in the original $(q,p)$
variables. We leave the details to the interested reader.

Finally we note that the core Toda flow can also be solved explicitly for any time $t$ in the bidiagonal formalism of \cite{LeSaTo2008}, so that in evaluating $$ \lim_{t \to\infty}\hat U_{-t} \circ U_t(q_0,p_0)$$ 
we could just as easily have worked in bidiagonal coordinates.

\bs
\noi{\bf Step 5.}

\ms

We assert that $W = (\hat W^\#)^{-1} \circ W^\#$ is $C^1(\RR^{2N})$. Indeed, $\hat W^\#$ is a diffeomorphism from $\RR^{2N}$ onto $\hat X^\#$, as can be seen directly from (\ref{35.2}) and (\ref{35.3}) and the fact that $\Phi$ is a diffeomorphism, coupled with the fact that
\[ (q_0, p_0) \mapsto \left( a_j \e - p_{j,0}/2,  b_j \e \frac{1}{2} e^{(q_{j,0} - q_{j+1,0})/2} , 2\le j\le N-1,  \sum_{i=2}^{N-1} q_{i,0} \right),\]
is a diffeomorphism. At the analytical level, the fact that $\hat W^\# \e \lim_{t \to \infty} U_{-t}^\# \circ \hat U_t$ is $C^1$ follows, alternatively, from the fact that $(q_0(t), p_0(t)) \e U_{-t}^\# \circ \hat U_t (q_0,p_0)$ is a $C^1$ function of $(q_0, p_0)$ for any finite $t$ by standard ODE methods, and then noting from (\ref{7.2}) that the derivative of $(q_0(t), p_0(t))$ with respect to $(q_0, p_0)$ is a linear system with exponentially decaying coefficients. But the same is true for (\ref{14.2}), and so $W^\# \e \lim_{t \to \infty} U_{-t}^\# \circ U_t$ is also $C^1$, and hence $W = (\hat W^\#)^{-1} \circ W^\#$ is $C^1$. We leave the details to the interested reader.

For the core Toda system $ (q,p)= (q_2, \ldots, q_{N-1}, p_2, \ldots, p_{N-2})$, let $\lambda_2 > \ldots > \lambda_{N-1}$ be the eigenvalues and $u_2(2), \ldots, u_{N-1}(2)$ be the first components of the normalized eigenvectors of the associated Toda matrix $L_F(q,p)$. Let $Z(\lambda) \e \det (L_F(q,  p) - \lambda)$. Then, as shown in \cite{DeLiNaTo1986}, and also, more directly in \cite{DeDuTr2021} (action-angle variables for the Toda flow are also given in \cite{KrVa2002}),
\beqn
\begin{array}{llll}
\theta_k \e \ln \left( \frac{u_k(2)}{u_{N-1}(2)} \ \Big| \frac{Z'(\lambda_k)}  {Z'(\lambda_{N-1})} \Big|^{1/2}\right) \ ,\quad 2 \le k \le N-2 \\
\theta_{N-1} \e \frac{1}{N-2} (  q_2 + \ldots +  q_{N-1}) ,\\
\tilde \lambda_k \e \lambda_k - \frac{1}{N-2} \sum_{j=2}^{N-1} \lambda_j \ , \quad 2 \le k \le N-2  ,\\
\tilde \lambda_{N-1} \e  p_2 + \ldots +  p_{N-1}
\end{array}
\eeqn
are action-angle variables for the core Toda flow,
\beql{43.2}
\{ \theta_i, \theta_j \} \e 0 \ , \{ \tilde \lambda_i , \tilde \lambda_j \} \e 0, \{\theta_i , \tilde \lambda_j \} = \delta_{ij} \ , 2 \le i, j \le N-1 \ .
\eeq
On the other hand, for additional variables $q_1, p_1,q_N,p_N$, if
\beqn
\begin{array}{ll}
\theta_1 \e - \frac{1}{2}  p_1 \ , \quad  \theta_N = \frac{1}{c}  p_N, \\
\tilde \lambda_1 \e \frac{1}{2}  p^2 + c q_1 \ , \quad \tilde \lambda_N \e \frac{1}{2}  p_N^2 - c  q_N ,
\end{array}
\eeqn
it follows that (\ref{43.2}) holds for $1 \le i, j \le N$.  Also from (\ref{31.7}),
\beql{43.4}
H_c^d \e \tilde \lambda_1 + \tilde \lambda_n + 2 \sum_{i=2}^{N-1} \lambda_k^2 \ .
\eeq
Now, as $\tilde \lambda_{N-1} = - 2( a_2 + \ldots + a_{N-1}) = - 2 \tr L_F \e - 2(\lambda_2 + \ldots + \lambda_{N-1})$, we see that $\lambda_k \e \tilde \lambda_k - \frac{1}{2(N-2)} \tilde \lambda_{N-1} , 2 \le k \le N-2$, and then solving for $\lambda_{N-1}$, we find
\[ \lambda_{N-1} \e - \frac{1}{2(N-2)} \tilde \lambda_{N-1} - \sum_{j=2}^{N-2} \tilde \lambda_i\ . \]
Substitution into (\ref{43.4}) gives
\[ H_c^d \e \tilde \lambda_1 + \tilde \lambda_N + 2 \sum_{k=2}^{N-2} (\tilde \lambda_k)^2 + 2\left( \sum_{k=2}^{N-2} \tilde \lambda_k\right)^2 + \frac{1}{2(N-2)^2} \tilde \lambda_{N-1}^2 \]
for which
\beql{44.1}
\begin{array}{lll}
\{ \theta_i , H_c^d \} \e 1 \ , \quad i = 1 \ \hbox{or} \ N \\
\{ \theta_i , H_c^d \} \e 4 \left( \tilde \lambda_i + \sum_{j=2}^{N-2} \tilde \lambda_j \right) \ , \quad 2 \le i \le N-2, \\
\{ \theta_{N-1}, H_c^d \} \e \frac{1}{(N-2)^2} \tilde \lambda_{N-1}
\end{array}
\eeq
so that the $\theta_i$'s move linearly under $H_c^d$. Also clearly,
\beql{44.4}
 \{ \tilde \lambda_i , H_c^d \} \e 0 \ , \quad 1 \le i \le N \ .
\eeq
Now for the flow $\hat U_t(\hat q_0, \hat p_0) \e (\hat q(t) , \hat p(t))$ ,
\[ \frac{d}{dt} f(\hat q(t), \hat p(t)) \e \{ f , H_c^q \}(\hat q(t) , \hat p(t))\]
for any $f: \RR^{2N} \to \RR$ and it follows from (\ref{44.1}) and (\ref{44.4}) that $\{\theta_i\}_{i=1}^N, \{ \tilde \lambda_i\}_{i=1}^N$ are action-angle variables for $H_c^d$.

Set
\beqn
\Theta_i \e \theta_i \circ W \ , \quad \Lambda_i \e \tilde \lambda_i \circ W \ , \quad 1 \le i \le N. \\
\eeqn
Then as in (\ref{12.1}), $\Theta_i$ and $\Lambda_i$ are canonically conjugate. Furthermore, under the flow $U_t$ generated by $H_c$, for any $f : \RR^{2N} \to \RR$ and  $F \e f \circ W$,
\begin{equation}
\begin{aligned}
 \frac{d}{dt} F(U_t(q_0,p_0)) \e & \frac{d}{dt} ( f \circ W \circ U_t(q_0,p_0))\\
= &  \frac{d}{dt} ( f \circ \hat U_t (\hat q_0, \hat p_0) ) \ , \ (\hat q_0, \hat p_0) \e W(q_0, p_0)\\
= & \frac{d}{dt} f(\hat q(t), \hat p(t)) \ . 
\end{aligned}
\end{equation}
Hence under $U_t$,
\begin{equation}
\begin{aligned}
& \frac{d}{dt} \Theta_i (q_0,p_0) \e \frac{d}{dt} \theta_i(W(q_0, p_0)) \e 1 \ , \quad i=1,N, \\
& \frac{d}{dt} \Theta_i (q_0,p_0) \e \frac{d}{dt} \theta_i(W(q_0, p_0))  \e 4\left(\Lambda_i(q_0,p_0) + \sum_{j=2}^{N-2} \Lambda_j(q_0,p_0)\right),\,\, \quad 2 \le i \le N-2, \\
& \frac{d}{dt} \Theta_{N-1} (q_0,p_0) \e \frac{1}{(N-2)^2} \Lambda_{N-1} (q_0, p_0) \ .  \\
\end{aligned}
\end{equation}
Also,
\beqn
\frac{d}{dt} \Lambda_i(q_0, p_0) \e 0 \ , \quad 1 \le i \le N \ .
\eeqn
Thus $\{ \Theta_i\}_{i=1}^N$, $\{\Lambda_i\}_{i=1}^N$ are action-angle variables for $H_c$.

Note that by a general and simple argument, the relations $\{ \Theta_i, \Lambda_j \} \e \delta_{ij}$ imply that $\Lambda_1, \ldots, \Lambda_n$ are 
functionally independent, which is equivalent to the statement that their gradients are linearly independent on an open, dense set. Thus $H_c$ is integrable in the sense of Liouville.

Finally note that every Hamiltonian of the form
\[ H \e \frac{1}{2} p^2 + V(q) \ \ (q,p) \in \RR^2\]
generates a flow that can be expressed in Lax-pair form. Indeed, for
\beql{45.3}\hat L \e \begin{pmatrix} p & 2 V(q) \\ 1 & - p \end{pmatrix} \ , \quad \hat B \e \begin{pmatrix} 0 & V'(q) \\ 0 & 0 \end{pmatrix},
\eeq
a simple computation shows that
\beql{45.4}
\dot{q} \e p \ , \quad \dot{p} \e - V'(q) \quad \Leftrightarrow \quad \frac{d}{dt} \hat L \e [ \hat L , \hat B ] \ .
\eeq
Let $(\hat q_k(t) , \hat p_k(t), 1 \le k \le N$, solve the flow generated by $H_c^d$ as above. Set
\beqn
Q_k \e \hat q_k \circ W \ , \quad P_k \e \hat p_k \circ W \ , \quad 1 \le k \le N \ .
\eeqn
Then, for $2 \le k \le N-1$,
\[ \frac{d}{dt} Q_k \e \hat p _ k \circ W \e P_k ,\]
\[ \frac{d}{dt} P_k \e \big(e^{\hat q_{k-1} - \hat q_k} - e^{\hat q_k - \hat q_{k+1}} \big) \circ W\]
\[ \e e^{ Q_{k-1} - Q_k} - e^{Q_k -  Q_{k+1}}, \, \]
where $e^{Q_1 - Q_2} = e^{Q_{N-1}  - Q_n} = 0 $. Also
\[ \frac{d}{dt} Q_1 \e P_1 \ , \quad  \frac{d}{dt} P_1 \e - c\]
\[  \frac{d}{dt} Q_N \e P_N \ , \quad  \frac{d}{dt} P_N \e c \ . \]
Finally, set
\beqn
{\mathcal L} \e \begin{pmatrix}
A_1 & B_2 & & & & & & \\
B_2 & \ddots &\ddots & & &0 & & \\
    & \ddots &  & B_{N-2} & & & & \\
    & & B_{N-2} & A_{N-1} & 0 & & & \\
     & & & 0 & P_1 & 2cQ_1 & & \\
    & & & & 1 & - P_1 & 0 & 0 \\
     & &0 & & & & P_N & -2c Q_N \\
     & & & & & & 1 & -P_N \\ \end{pmatrix} ,
\eeqn
and
\beqn
{\cal B} \e \begin{pmatrix}
0 & -B_2 & & & & & & \\
B_2 & \ddots &\ddots & & &0 & & \\
    & \ddots &  & -B_{N-2} & & & & \\
    & & B_{N-2} & 0 & 0 & & & \\
     & & & 0 & 0 & c & & \\
    & & & & 0 &0 & 0 & \\
     &0 & & & & & 0 & -c \\
     & & & & & & 0 & 0 \\ \end{pmatrix} ,
\eeqn
where
\beqn
\begin{array}{ll}
A_k \e - \frac{1}{2} P_k \ , \quad 2 \le k \le N-1, \\
B_k \e \frac{1}{2} e^{(Q_k - Q_{k+1})/2} \ , \quad 2 \le k \le N-2 \ .
\end{array}
\eeqn
Then

\centerline{  $(q_k(t) , p_k(t))_{k=1}^N$\ \  solve (\ref{14.2})}
\[ \Leftrightarrow \]
\[ \frac{d}{dt} {\cal L} \e [ {\cal L} , {\cal B} ] \ . \]
Thus the Hamilton equations for $H_c$ have a Lax-pair form.

\bs
\noi{\bf Remark:} Instead of using (\ref{34.4}), we could use (\ref{34.5}),
$U_t^\# \circ W^\# \e W^\# \circ U_t$ to display (\ref{14.2}) as a Lax-pair in another form. But now the analog ${\mathcal L}^\#$ and ${\mathcal B}^\#$  of ${\mathcal L}$ and ${\mathcal B}$ convey little information,
\[ {\mathcal L}^\# \e \begin{pmatrix}
P_1 & 0& & & & & & \\
0 & P_2 &\ddots & & &0 & & \\
    & \ddots & \ddots & 0 & & & & \\
    & & 0 & P_{N-1} & 0 & & & \\
     & & & 0 & P_1 & 2cQ_1 & & \\
    &0 & & & 1 & - P_1 & 0 & \\
     & & & & & & P_N & -2c Q_n \\
     & & & & & & 1 & -P_N \\ \end{pmatrix} \]
     and
\[ {\cal B}^\# \e \begin{pmatrix}
0 &0 & & & & & & \\
0 & \ddots &\ddots & & &0 & & \\
    & \ddots &  &0 & & & & \\
    & & 0 & 0 & 0 & & & \\
     & & & 0 & 0 & c & & \\
    & & & & 0 &0 &  & \\
     &0 & & & & & 0 & -c \\
     & & & & & & 0 & 0 \\ \end{pmatrix} \]
     where
\beql{48.1}
Q_k \e q_k^\# \circ W^\# \ , \quad P_k \e p_k^\# \circ W^\# \ , \quad 1 \le k \le N \ . \eeq
Here $(q^\# (t), p^\#(t)) = U^\#_t(q^\#_0,p^\#_0)$.

\vskip .1in
\textbf{Acknowledgments.}  Luen-Chau Li acknowledges the support from the Simons Foundation through grant \#585813. Tom Trogdon acknowledges the support of NSF grant DMS-1945652 and Carlos Tomei from CNPq and FAPERJ, Brazil. The authors would also like to thank Barry Simon for some helpful information.

\end{document}